\documentclass[9pt]{article} % For LaTeX2e
\usepackage{hyperref}
\usepackage{url}
\usepackage{smile}
\usepackage{graphicx} % more modern
\usepackage{algorithm}
\usepackage{algorithmic}
\usepackage{epstopdf}
\usepackage[margin=1.0in]{geometry}
\usepackage[export]{adjustbox}% http://ctan.org/pkg/adjustbox
\usepackage{mathtools, cuted}
\usepackage{bbm}
\usepackage{wrapfig}
\usepackage{subcaption}
\usepackage{caption}
\usepackage{enumitem}
\usepackage{tabularx}
\usepackage{smile}
\usepackage{tcolorbox}
\usepackage{enumitem}
\usepackage{mathtools, cuted}
\usepackage{caption}
\usepackage{subcaption}

\numberwithin{equation}{section}

\usepackage[nocompress]{cite}

\hypersetup{
    colorlinks=true,
    linkcolor=blue,
    anchorcolor=blue, 
    citecolor=blue,
}

\allowdisplaybreaks

\renewcommand{\frac}[2]{\tfrac{#1}{#2}}

\title{
First-order Policy Optimization for  Robust Policy Evaluation
    }
\author{
    Yan Li   \thanks{H. Milton Stewart School of Industrial and Systems Engineering, Georgia Institute of Technology, Atlanta, GA, 30332. (E-mail: \url{yli939@gatech.edu}).}
        \and 
    Guanghui Lan \thanks{H. Milton Stewart School of Industrial and Systems Engineering, Georgia Institute of Technology, Atlanta, GA, 30332. (E-mail: \url{george.lan@isye.gatech.edu}).}
}
\date{\vspace{-6ex}}

\begin{document}
{
\makeatletter
\addtocounter{footnote}{1} % to get dagger instead of star
\renewcommand\thefootnote{\@fnsymbol\c@footnote}%
\makeatother
\maketitle
}

\maketitle

\begin{abstract}
We adopt a policy optimization viewpoint towards policy evaluation for robust Markov decision process with $\mathrm{s}$-rectangular ambiguity sets. The developed method, named first-order policy evaluation (FRPE), provides the first unified framework for robust policy evaluation in both deterministic (offline) and stochastic (online) settings, with either tabular representation or generic function approximation. In particular, we establish linear convergence in the deterministic setting, and $\tilde{\cO}(1/\epsilon^2)$ sample complexity in the stochastic setting. FRPE also extends naturally to evaluating the robust state-action value function with $(\mathrm{s}, \mathrm{a})$-rectangular ambiguity sets. We discuss the application of the developed results for stochastic policy optimization of large-scale robust MDPs. 
\end{abstract}

%!TEX root = ./robust_pe.tex

\section{Introduction}

A robust Markov decision process (MDP)  $\cM_\cP \coloneqq \cbr{\cM_{\PP} : \PP \in \cP}$ consists of a collection of MDPs.
Each MDP $\cM_{\PP} = (\cS, \cA, c, \PP, \gamma)$ is a quintuple, where  
  $\cS$ and $\cA$ denote the finite state and action space, respectively;  
 $c: \cS \times \cA \to \RR$ denotes the cost function, which we assume without loss of generality that $c(s,a) \in [0,1]$ for all $(s,a)$; $\gamma$ denotes the discount factor;
and $\PP \in \Delta_{\cS}^{\abs{\cS}  \abs{\cA}}$ denotes the transition kernel,
where $\Delta_{\cS}$ corresponds to the probability simplex over $\cS$.
Notably MDPs within $\cM_{\cP}$ differ from each other only in their respective transition kernels. 
The set of potential transitional kernels $\cP$ is referred to as the ambiguity set.
%We adopt the following notation going forward. 

For $\PP \in \Delta_{\cS}^{\abs{\cS} \abs{\cA}}$, we use 
$\PP_s \in \Delta_{\cS}^{\abs{\cA}}$ to denote its component corresponding to $s \in \cS$, 
and $\PP_{s,a} \in \Delta_{\cS}$ to denote its component of $(s,a) \in \cS \times \cA$.
Similarly, for $\DD \in \Delta_{\cS}^{\abs{\cA}}$,  $\DD_{a} \in \Delta_{\cS}$ denotes its component of $a \in \cA$.
%In addition,  $\PP \in (\Delta_{\cS})^{\abs{\cS} \times \abs{\cA}}$ denotes the transition kernel, with $\Delta_{\cS}$ being the probability simplex over $\cS$. 
%Accordingly, $\PP_{s,a} \in \Delta_{\cS}$ corresponds to the conditional distribution of the next state upon making an action $a \in \cA$ at state $s \in \cS$,
%and we let $\PP_s = [\PP_{s,a_1}, \ldots, \PP_{s, a_{\abs{\cA}}}]$.
%The set of MDPs differ from each other only in their respective transition kernels. 
Then the standard value function $V^{\vartheta}_{\PP}: \cS \to \RR$ of a randomized stationary policy $\vartheta$ with respect to $\cM_{\PP}$  is defined as 
\begin{align}\label{eq_standard_value_function}
V^{\vartheta}_{\PP}(s) \coloneqq \EE \sbr{\tsum_{t=0}^\infty \gamma^t c(s_t, a_t) \big| s_0 = s, a_t \sim \vartheta(\cdot|s_t), s_{t+1} \sim \PP_{s_t, a_t} }, ~~ \forall s \in \cS.
\end{align}

We consider policy evaluation for robust Markov decision process with $\mathrm{s}$-rectangular ambiguity sets \cite{wiesemann2013robust}.
That is, for a given policy $\vartheta$, we aim to find worst-case value function  defined as
\begin{align}\label{eq_def_robust_value}
\textstyle
V^{\vartheta}_r (s) \coloneqq \max_{\PP \in \cP} V^{\vartheta}_{\PP}(s), ~ \forall s \in \cS,
\end{align}
where the set of potential transition kernels $\cP \in \Delta_{\cS}^{\abs{\cS}  \abs{\cA}}$ decomposes as 
\begin{align}\label{eq_def_s_rectangular}
\cP = \times_{s \in \cS}  \cP_{s},  ~ \text{where}~ \cP_{s} \subseteq \Delta_{\cS}^{\abs{\cA}}.
\end{align}
%Ambiguity sets with above decomposition are first introduce in \cite{wiesemann2013robust} and referred to as the $\mathrm{s}$-rectangular set.
%We refer as the marginal ambiguity set in the ensuing discussion.
We assume $\cP_{s}$ takes the form of 
\begin{align}\label{def_ambiguity_set_structure}
\cP_{s} = \cbr{
(1-\zeta) \overline{\PP}_{s} + \zeta \DD: ~ \DD \in \cD_{s}
},
\end{align}
where $\overline{\PP}$ is a fixed but possibly unknown transition kernel,   $\cD_{s}$ is  a pre-specified convex subset of $\Delta_{\cS}^{\abs{\cA}}$, and $\zeta \in [0,1]$ is a parameter of user choice. 
One can also specify $\zeta$ in a potentially state-dependent manner.
%We also denote 
%$\cD = \times_{s \in \cS} \cD_{s}$.
%\begin{align}
%\label{eq_cD_set}
%\cD = \times_{s \in \cS} \cD_{s}.
%\end{align}
We next discuss two examples of robust MDPs that can be modeled by \eqref{def_ambiguity_set_structure}.
%We do not assume any information on $\PP$ except that we can sample from it.

%\begin{example}[Robust MDP with Offline Data]\label{example_offline_rmdp}
\begin{example}[Offline Robust MDP]\label{example_offline_rmdp}
One of the most classical applications of robust MDP deals with planning with a historical dataset, collected by prior interactions with an environment of unknown transition kernel $\PP^{\mathrm{N}}$.
The goal is to perform robust planning in response to the inherent uncertainty associated with data collection. 
In this setting, one would first construct for each state a set of plausible models $\cD_{s}$, such that  $ \PP^{\mathrm{N}}_{s} \in \cD_{s}$ with sufficient confidence.
Then the planner solves the robust MDP 
\begin{align*}
\textstyle
\min_{\vartheta} V^{\vartheta}_r(s), ~ \forall s \in \cS,
\end{align*}
where the ambiguity set is defined by \eqref{eq_def_s_rectangular}, with $\zeta = 1$ in \eqref{def_ambiguity_set_structure}, and $\overline{\PP}$ being an arbitrary transition kernel.
% and solves the robust MDP problem with ambiguity set $\cP$. 
Such an approach can be found in early studies of robust MDPs \cite{nilim2005robust, iyengar2005robust, wiesemann2013robust}.
\end{example}

%\begin{example}[Robust MDP without Offline Data]
\begin{example}[Online/Hybrid Robust MDP]\label{online_robust_mdp}
Another emerging class of applications of robust MDP seeks robustness against environment shift \cite{roy2017reinforcement, liu2022distributionally, wang2022policy}.
Namely, the goal is to learn a robust policy that performs well in both training and potentially altered testing environment. 
Notably both training and testing environment can be unknown, except that one can access samples from the training environment, and the testing environment should be close to the training one. 
A typical example of this scenario is the so-called sim-to-real process, where the target is to learn a robust policy within a simulated environment that behaves reasonably well in reality despite the modeling errors associated with simulation.
It can be seen that \eqref{def_ambiguity_set_structure} creates a mechanism of robustness after an immediate reformulation
\begin{align*}
\cP_{s} = \cbr{
 \overline{\PP}_{s} + \zeta (\DD - \overline{\PP}_{s}): ~ \DD \in \cD_{s}
}.
\end{align*}
Here $\overline{\PP}$ denotes the transition kernel of the training environment, while $\cU_\zeta = \cbr{\zeta (\DD - \overline{\PP}_{s}): ~ \DD \in \cD_{s}}$ corresponds to the potential differences between testing and training environment.
Clearly $\zeta$ provides a way to adjust the size of the ambiguity set (and consequently how much environment change one can tolerate) based on our robustness preference.
The flexibility of choosing $\cD_{s}$ also provides a convenient way to utilize partial offline data or domain knowledge.
Namely, for a particular state, if offline data or domain knowledge is available, one can construct $\cD_{s}$ by exploiting this knowledge.
Otherwise one may set $\cD_{s}$ as $\Delta_{\cS}^{\abs{\cA}}$. 
\end{example}

\begin{remark}
Our discussions focus on the robust value function with $\mathrm{s}$-rectangular sets. 
Yet it can be readily seen that the developed methods in this manuscript extend naturally to evaluating the robust Q-function with $(\mathrm{s}, \mathrm{a})$-rectangular sets.
Note that for $\mathrm{s}$-rectangular sets, there is no natural notion of robust Q-function. 
%, and the value function for $\mathrm{s}$-rectangular sets.
%For the latter class of ambiguity sets, there is no natural notion of robust Q-function.
\end{remark}

Going forward we refer to the setting where $\overline{\PP}$ is known as the deterministic setting, and 
where it is unknown as the stochastic setting.
Clearly, offline robust MDP can be categorized into deterministic setting while online robust MDP belongs to the stochastic setting.

The problem of robust policy evaluation \eqref{eq_def_robust_value} serves as the essential step of policy optimization for robust MDPs. 
Both classical methods of robust policy iteration and recently developed first-order policy optimization (i.e., policy gradient) methods \cite{wang2022policy, li2022first} hinge upon an accurate estimation of the robust value function. 
For deterministic setting with known ambiguity sets, 
robust policy evaluation can be efficiently done for small-sized problems. 
In particular,  recursive applications of the robust Bellman operator  yields linear convergence to the robust value function \cite{nilim2005robust, iyengar2005robust} as the robust Bellman operator is a contraction.
%much similar to how non-robust Bellman operator can be used for evaluating standard value function \eqref{eq_standard_value_function}.
% Both methods  depend critically on the contraction property of the Bellman operators. 
This property is further exploited in the stochastic setting, 
where an $\cO(1/\epsilon^2)$ sample complexity is obtained with a robust temporal difference learning method \cite{li2022first}.
Even for this simple tabular problem, it remains unclear whether robust policy evaluation requires additional samples compared to its non-robust counterpart.

For large-scale problems,  robust policy evaluation becomes challenging even for the deterministic setting.
For computational consideration one needs to introduce certain level of approximation to handle large state space.
One prior approach focuses on linear function approximation \cite{tamar2014scaling}, where the robust Bellman operator is replaced by its robust projected counterpart  applied to the linearly-parameterized robust value estimation. 
Unfortunately, this modification destroys the contraction property of the resulting fixed point update and leads to divergence.
Indeed a rather unrealistic assumption on the discount factor and the ambiguity set is necessary to retain  contraction  \cite{tamar2014scaling}.
% Yet modification destroys the contraction properties of the resulting fixed point update and leads to divergence \cite{tamar2014scaling}. 
% A necessary and sufficient condition is introduced in \cite{tamar2014scaling}  to recover contraction. Unfortunately, this condition places unrealistic requirement on the discount factor and the ambiguity set, and is generally not verifiable. 
Another natural approach seeks to directly minimize Bellman residual. 
This leads to a linear system for standard policy evaluation problem \eqref{eq_standard_value_function} but becomes easily non-convex for robust policy evaluation \eqref{eq_def_robust_value},  due to the non-linearity of the robust Bellman operator \cite{roy2017reinforcement}. 
For the stochastic setting, there seems to be no evaluation method without a similar assumption \cite{roy2017reinforcement}.
%It appears that there is no stable and globally convergent method of robust policy evaluation for large-scale problems. 

An alternative approach, yet largely dismissed, is to consider robust policy evaluation as a policy optimization problem in its own right.
%The connection between robust policy evaluation and policy optimization appears to be known.
The connection of these two problems appears in various prior development.
It is discussed in \cite{ruszczynski2010risk, shapiro2021distributionally} that certain nested formulations of robust MDP can be viewed as a stochastic game between the planner and the nature, with the nature's optimal policy defining the robust value of the planner. 
From our perspective the omission of this viewpoint may be attributed to two reasons.
First,  early applications of offline robust MDPs (Example \ref{example_offline_rmdp}) can already be solved for small-sized problems. 
This is further compounded by the fact that in the minimax game the nature has a continuous action space, and parameterizing the nature's policy seems to already require a memory of the same size as the transition~kernel.

In this manuscript we develop the first unified and globally convergent framework of robust policy evaluation for both deterministic and stochastic settings, with  either tabular representation or generic function approximation. 
Enabled by this framework we establish a set of computational results that appear to be completely new for robust policy evaluation problems, and greatly enhance the applicability of robust MDPs to large-scale problems. 
Our main contributions are summarized into the following aspects.

%In this manuscript we develop the first stable and unified framework of robust policy evaluation for both deterministic and stochastic settings, with  either tabular representation or generic function approximation. 
%Enabled by this framework we establish a set of computational results that appear to be completely new for robust policy evaluation problems, and greatly enhance the applicability of robust MDPs to large-scale problems. 
%Notably we adopt a policy optimization viewpoint towards robust policy evaluation, and consequently our method development does not require the contraction property.
%The connection between robust policy evaluation and policy optimization appears to be known.
%It is shown in \cite{ruszczynski2010risk, shapiro2021distributionally} that certain nested formulations of robust MDP can be viewed as a minimax stochastic game between the planner and the nature, with the nature's optimal policy defining the robust value of the planner. 
%Nevertheless, this viewpoint has largely been dismissed for algorithmic development.
%From our perspective this might be attributed to two reasons.
%First,  early applications of robust MDPs (Example \ref{example_offline_rmdp}) can already be solved for small-sized problems by using exact fixed point iterations. 
%This is further compounded by the fact that in the minimax  game the nature has a continuous action space, and parameterizing the nature's policy seems to already require a memory of the same size as the transition kernel. 
%Our main contributions are summarized into the following aspects. 

First, we formulate the robust policy evaluation problem \eqref{eq_def_robust_value} into a Markov decision process of the nature, whose action space corresponds to the ambiguity set of the robust MDP. 
We then develop a method named first-order robust policy evaluation (FRPE), which only requires solving a standard policy evaluation problem \eqref{eq_standard_value_function} to improve the policy of the nature.
In a nutshell, FRPE operates much similarly to the dual averaging method \cite{nesterov2009primal} and its recent adaptation to policy optimization for general state and action spaces \cite{lan2022policy}. Accordingly, despite being a policy based method, FRPE can be implemented in a fully value-based manner and avoid explicit policy parameterization. 
%Consequently, no policy parameterization is required for FRPE.

In the deterministic setting, we show that FRPE converges linearly to the robust value function for tabular problems with small state spaces. 
Each step of FRPE only requires a simple subroutine for standard policy evaluation with a known transition kernel.
%Notably, the iteration complexity is independent of the size of the action space for the planner, despite the nature's MDP having a continuous action space.
We then extend FRPE to settings where function approximation is employed to handle large state space. 
In this case FRPE converges at the same linear rate, up to the approximation error produced by the standard evaluation subroutine. 
Consequently, FRPE substantially improves existing methods of robust policy evaluation for large-scale problems as it does not require any restrictive assumption, and can be combined with general function approximation schemes.
%as no assumption is required for convergence, and is applicable to general function approximation schemes. 

In the stochastic setting, we develop a stochastic variant named SFRPE.
For tabular problems, we establish an $\cO\rbr{\frac{\zeta^2 \abs{\cS}}{(1-\gamma)^5 \epsilon^2} + \frac{\abs{\cS}}{1-\gamma} \log (\frac{1}{\epsilon})}$ sample complexity for SFRPE to return an estimated robust value function with its bias bounded by $\epsilon$.
%is $\epsilon$-accurate in expectation (Definition \ref{def_acc_certificate}). 
Notably, the established sample complexity is independent of the continuous action space of the nature. 
We then establish an $\cO\rbr{\frac{\zeta^2 \abs{\cS}}{(1-\gamma)^5 \epsilon^2} + \frac{\abs{\cS}}{(1-\gamma)^3} \log (\frac{1}{\epsilon})}$ sample complexity for 
returning $\epsilon$-accurate estimator with high probability.
As a consequence the robust policy evaluation does not require additional samples compared to standard policy evaluation, provided the size $\zeta$ of the ambiguity set is bounded by $\cO(1-\gamma)$. 
For large-scale problems, we further incorporate SFRPE with linear function approximation and establish a similar $\cO\rbr{{\zeta^2}/{\epsilon^2} + \log({1}/{\epsilon})}$ (resp. $\cO\rbr{{\zeta^2}/{\epsilon^2} + {1}/{\epsilon^2}}$) sample complexity for returning an $\epsilon$-accurate estimator in expectation (resp. in high probability).
The dependence on the size of the ambiguity set clearly delineates the price of robustness for robust policy evaluation. 
In particular, setting $\zeta =0$ recovers a tight sample complexity for standard policy evaluation. 
To the best of our knowledge, all the obtained sample complexities appear to be new.
Importantly, with linear function approximation we have established the first known sample complexity of robust policy evaluation beyond the tabular setting.
SFRPE is flexible enough and can be combined even with nonlinear function approximation, provided certain off-policy evaluation problem can be solved to a prescribed target precision.

Finally, we demonstrate a direct application of SFRPE for solving large-scale robust MDPs.
 Combining a recently developed stochastic policy optimization method for robust MDPs \cite{li2022first} and SFRPE, we establish an $\tilde{\cO}(1/\epsilon^2)$ sample complexity for both tabular problems, and for large-scale problems with log-linear policy class. 
Notably the latter sample complexity has not been reported for this problem class.

The rest of the manuscript is organized as follows.
Section \ref{sec_pe_as_po}  formulates the MDP for robust policy evaluation and studies a few of its useful structural properties. 
Section \ref{sec_deterministic} then introduces FRPE for the deterministic setting and establishes its iteration complexities.
Section \ref{sec_stochastic} develops the stochastic FRPE and its sample complexities. 
 Finally, we conclude in section \ref{sec_conclusion}, and discuss a few future directions in robust policy evaluation and optimization  enabled by FRPE.

%\yan{Key points to make:
%
%
%
%\begin{itemize}
%\item need to mention one of the biggest strength --
%one does not need to store the nature's policy -- all we need is store weighted value and function approximation can be used -- making it truly scalable
%\item the first term of the sample complexity corresponds to the price we pay for robustness -- notably when $\zeta = \cO(1-\gamma)$ there is no price to pay for robustness!
%\end{itemize}
%
%}

%!TEX root = ./robust_pe.tex
%\newpage

\section{Robust Policy Evaluation as Policy Optimization}\label{sec_pe_as_po}
This section adopts a policy optimization viewpoint towards policy evaluation, by first formulating the robust policy evaluation problem as a Markov decision problem of the nature.
We then identify a few key structures of the formulated MDP that will prove useful for our subsequent development. 

Consider a MDP of nature, denoted by $\mathfrak{M}$,  defined as follows.
The state space is given by $\cS$, and
the set of possible actions  at any state $s \in \cS$ is given by $ \cD_s \subseteq \Delta_{\cS}^{\abs{\cA}}$ (cf. \eqref{def_ambiguity_set_structure}).
%We will write $\DD = [\DD_{a_1}, \ldots, \DD_{a_{\abs{\cA}}}]$ for any $\DD \in \Delta_{\cS}^{\abs{\cA}}$.
%Equivalently, any possible action $\DD \in \cD_s$ specifies $\abs{\cA}$ elements in $\Delta_{\cS}$,  each denoted as $\DD_{a}$ for every $a \in \cA$.
At state $s \in \cS$, upon making an action $\DD \in \cD_s$, the conditional distribution of the next state $s' \in \cS$ is given by 
\begin{align}\label{transit_kernel_of_nature_mdp}
\mathfrak{P}(s' | s, \DD) \coloneqq  \tsum_{a \in \cA} \sbr{(1-\zeta) \overline{\PP}_{s,a}(s')  + \zeta \DD_{a}(s')} \vartheta(a|s).
\end{align}
%Clearly, the above transition kernel $\mathfrak{P}$ is affine with respect to the action of the nature  $\DD$. 
Finally, the cost function associated with $(s, \DD)$ for any $\DD \in \cD_s$ is given by 
\begin{align*}
\mathfrak{C}(s) \coloneqq - \tsum_{a \in \cA} \vartheta(a|s) c(s,a),
\end{align*}
 and the discount factor is set as $\gamma$.

A non-randomized policy of the nature is denoted as $\pi: \cS \to \cD_s$. 
%It is clear that $\pi$ uniquely determines $\DD^{\pi} \in \cD$ defined in \eqref{eq_cD_set}.
%Let us denote $\DD^{\pi(s)} = \pi(s)$, and 
%Let us define $\DD^{\pi(s)}  = \pi(s)$ for any policy $\pi$.
For any policy $\pi$ and $s \in \cS$, let us define $\DD^{\pi(s)} \equiv \pi(s) \in \cD_s$.
For notational clarity we will sometimes use these two quantities interchangeably. 
We then define
\begin{align}\label{kernel_defined_by_nature_policy}
\PP^{\pi}_{s,a} \coloneqq (1-\zeta) \overline{\PP}_{s,a}  + \zeta \DD^{\pi(s)}_{a}, ~ (s,a) \in \cS \times \cA.
\end{align}
%as the transition kernel of the original planner when the nature's policy is $\pi$.
Consequently, from \eqref{transit_kernel_of_nature_mdp} it holds that 
\begin{align}\label{state_transit_of_nature_given_policy}
\mathfrak{P}(s'|s,  \pi(s)) =  \tsum_{a \in \cA} \sbr{(1-\zeta) \overline{\PP}_{s,a}(s')  + \zeta \DD^{\pi(s)}_{a}(s')} \vartheta(a|s)
= \tsum_{a \in \cA} \PP^{\pi}_{s,a}(s') \vartheta(a|s).
\end{align}
We define the value function of policy $\pi$ as 
\begin{align*}
%\label{eq_def_value_func_nature}
\cV^{\pi} (s) \coloneqq 
\EE \sbr{\tsum_{t=0}^\infty \gamma^t \mathfrak{C}(s_t) \big| s_0 = s, s_{t+1} \sim \mathfrak{P}(\cdot| s_t, \pi(s_t) ) , t \geq 0}, ~~ \forall s \in \cS,
\end{align*}
and the goal of the nature is to find the optimal policy of 
\begin{align}\label{eq_def_optmal_value_nature}
\textstyle
\min_{\pi \in \Pi} \cV^{\pi} (s),
\end{align}
where $\Pi: s \mapsto \cD_s$ is the set of non-randomized stationary policies of the nature.

\begin{lemma}\label{lemma_value_correspondence}
For any $\pi \in \Pi$, we have 
\begin{align}\label{eq_nature_value_as_player_value}
\cV^{\pi}(s) = - V^{\vartheta}_{\PP^{\pi}}(s), ~ \forall s \in \cS,
\end{align}
with $V^{\vartheta}_{\PP^{\pi}}$ defined in \eqref{eq_standard_value_function}.
In addition, let $\cV^*$ denote the optimal value function of \eqref{eq_def_optmal_value_nature}.
Then 
\begin{align}\label{nature_opt_as_robust_value}
\cV^*(s) = - V^{\vartheta}_r(s), ~ \forall s \in \cS,
\end{align}
where $V^{\vartheta}_r$ is defined as in \eqref{eq_def_robust_value}.
\end{lemma}

\begin{proof}
It is clear that the $\cV^{\pi}$ satisfies the following dynamic programming equation 
\begin{align*}
\cV^{\pi} (s) & = \mathfrak{C}(s) + \gamma \tsum_{s' \in \cS} \mathfrak{P}(s'|s,  \pi(s)) \cV^{\pi}(s') \\
 & = - \tsum_{a \in \cA} \vartheta(a|s) c(s,a)
 + \gamma  \tsum_{s' \in \cS} \tsum_{a \in \cA} \PP^{\pi}_{s,a}(s') \vartheta(a|s) \cV^{\pi}(s'), ~ \forall s \in \cS.
\end{align*}
where the last equality follows from \eqref{state_transit_of_nature_given_policy}.
The above relation implies that $- \cV^{\pi}$ is the fixed point of operator 
\begin{align*}
(\cT^{\pi} V)(s) = 
\tsum_{a \in \cA} \vartheta(a|s) c(s,a)
 + \gamma  \tsum_{s' \in \cS} \tsum_{a \in \cA} \PP^{\pi}_{s,a}(s') \vartheta(a|s) V(s'), ~ \forall s \in \cS.
\end{align*}
On the other hand, it is known that $V^{\vartheta}_{\PP^{\pi}}$ is the unique fixed point of $\cT^{\pi}$, from which
we obtain \eqref{eq_nature_value_as_player_value}.
% That is, the value function of the nature's policy $\cV^{\pi}$ corresponds to the negative value function of the policy $\pi$ within $\cM_{\PP^{\pi}}$.
In addition,  Bellman optimality condition of MDP \eqref{eq_def_optmal_value_nature} yields  
\begin{align*}
\cV^*(s) & = \min_{\DD \in \cD_s} \mathfrak{C}(s)  + \gamma \tsum_{s' \in \cS} \mathfrak{P}(s' |s, \DD) \cV^*(s') \\ 
& = \min_{\DD \in \cD_s} -  \tsum_{a \in \cA} \vartheta(a|s) c(s,a) 
+ \gamma \tsum_{s' \in \cS} \tsum_{a \in \cA}\sbr{(1-\zeta) \overline{\PP}_{s,a}(s')  + \zeta \DD_{a}(s')}\vartheta(a|s) \cV^*(s') \\
& = \min_{\PP \in \cP_s} -  \tsum_{a \in \cA} \vartheta(a|s) c(s,a) 
+ \gamma \tsum_{s' \in \cS} \tsum_{a \in \cA} \PP_{a}(s') \vartheta(a|s) \cV^*(s') , ~ \forall s \in \cS.
%\\
%& =  - \tsum_{a \in \cA} \vartheta(a|s) c(s,a) 
%+ \gamma  \tsum_{a \in \cA}   \vartheta(a|s)  \tsum_{s' \in \cS} \min_{\PP_{s,a} \in \cP_{s,a}} \PP_{s,a}(s') \cV^*(s'), ~ \forall s \in \cS.
\end{align*}
Clearly, $-\cV^*$ is the fixed point of operator 
\begin{align}\label{def_robust_ballmen_op}
(\cT V)(s) = \max_{\PP \in \cP_s} \tsum_{a \in \cA} \vartheta(a|s) c(s,a) 
+ \gamma \tsum_{s' \in \cS} \tsum_{a \in \cA} \PP_{a}(s') \vartheta(a|s) V(s'), ~ \forall s \in \cS. 
\end{align}
On the other hand, it is well known that $V^{\vartheta}_{r}$ is the unique fixed point of $\cT V$ \cite{wiesemann2013robust}. 
Consequently we obtain \eqref{nature_opt_as_robust_value}.
%.
%That is, the optimal value function of the nature \eqref{eq_def_value_func_nature} corresponds to the negative robust value function $V^{\pi}_r$ of the policy,
%as both are the (unique) solution of the above dynamic programming equation.
\end{proof}

In view of the above observations, the robust policy evaluation problem \eqref{eq_def_robust_value} can be equivalently solved by solving a Markov decision process \eqref{eq_def_optmal_value_nature}  of nature with finite state space and continuous action space. 
To this end, let us define the following problem:
\begin{align}\label{policy_opt_obj_nature}
\textstyle
\min_{\pi \in \Pi} \cbr{f(\pi) \coloneqq \EE_{s \sim \rho} \sbr{\cV^{\pi}(s)}},
\end{align}
where $\rho$ is any distribution with full support over $\cS$.
Our end goal is to develop efficient methods that can be applied to solve \eqref{policy_opt_obj_nature}.

The state-action value function of the nature, also know as the Q-function, is defined by
\begin{align}\label{def_q_func_nature}
\cQ^{\pi}(s, \DD) & \coloneqq 
\EE \sbr{\tsum_{t=0}^\infty \gamma^t \mathfrak{C}(s_t) \big| s_0 = s, s_1 \sim   \mathfrak{P}(\cdot| s, \DD ), s_{t+1} \sim \mathfrak{P}(\cdot| s_t, \pi(s_t) ), t \geq 1 } . 
%\\
%& =  \mathfrak{C}(s) + \gamma \EE_{s' \sim   \mathfrak{P}(\cdot| s, \DD )} \sbr{\cV^{\pi}(s')} \\
%& = \mathfrak{C}(s) +  
% \gamma \tsum_{s' \in \cS} \tsum_{a \in \cA}\sbr{(1-\zeta) \overline{\PP}_{s,a}(s')  + \zeta \DD_{s,a}(s')}\vartheta(a|s) \cV^{\pi}(s')  \\
%& =  \mathfrak{C}(s) +  
% \gamma  \tsum_{a \in \cA} \vartheta(a|s) 
% \inner{(1-\zeta) \overline{\PP}_{s,a} + \zeta \DD_{s,a}}{\cV^{\pi}} \\
% & = 
% \mathfrak{C}(s) + \gamma \inner{(1-\zeta) \overline{\PP}_s + \zeta \DD}{\cV^{\pi}_{\vartheta, s}}, 
\end{align}
Clearly one also has 
\begin{align}\label{relation_q_and_v}
\cQ^\pi(s,\DD) = \mathfrak{C}(s) + \gamma \EE_{s' \sim \mathfrak{P}(\cdot | s,\DD)} \sbr{\cV^{\pi}(s')}.
\end{align}
%for any $ \DD \in \cD_s$, 
%where in the last equality we define $\cV^{\pi}_{\vartheta, s} = \vartheta(\cdot|s) \otimes \cV^{\pi} \in \RR^{\abs{\cA} \abs{\cS}}$.
We next show that $Q^{\pi}(s, \cdot)$ is indeed an affine function over $\cD_s$, an immediate yet important property that we will exploit in the ensuing development. 

\begin{proposition}\label{prop_q_structure}
For any $\DD \in \cD_s$, we have 
\begin{align*}
\cQ^{\pi}(s, \DD)
=  \mathfrak{C}(s) + \gamma \inner{(1-\zeta) \overline{\PP}_s + \zeta \DD}{\cV^{\pi}_{\vartheta, s}}, 
\end{align*}
where $\cV^{\pi}_{\vartheta, s} \coloneqq \vartheta(\cdot|s) \otimes \cV^{\pi} \in \RR^{ \abs{\cS} \abs{\cA}}$.
%and $\overline{\PP}_s \coloneqq [\overline{\PP}_{s, a_1}, \ldots, \overline{\PP}_{s, a_{\abs{\cA}}}] \in \Delta_{\cS}^{\abs{\cA}}$.
\end{proposition}

\begin{proof}
We have 
\begin{align*}
\cQ^{\pi}(s, \DD)
& =  \mathfrak{C}(s) + \gamma \EE_{s' \sim   \mathfrak{P}(\cdot| s, \DD )} \sbr{\cV^{\pi}(s')} \\
& = \mathfrak{C}(s) +  
 \gamma \tsum_{s' \in \cS} \tsum_{a \in \cA}\sbr{(1-\zeta) \overline{\PP}_{s,a}(s')  + \zeta \DD_{a}(s')}\vartheta(a|s) \cV^{\pi}(s')  \\
& =  \mathfrak{C}(s) +  
 \gamma  \tsum_{a \in \cA} \vartheta(a|s) 
 \inner{(1-\zeta) \overline{\PP}_{s,a} + \zeta \DD_{a}}{\cV^{\pi}} \\
 & = 
 \mathfrak{C}(s) + \gamma \inner{(1-\zeta) \overline{\PP}_s + \zeta \DD}{\cV^{\pi}_{\vartheta, s}},
\end{align*}
which completes the proof.
\end{proof}

%\yan{need to define $d_{\rho}^{\pi}$ within the perf diff lemma}
Our ensuing discussions repeatedly make use of the discounted visitation measure, defined as $d_{\rho}^{\pi}(s) = (1-\gamma) \tsum_{s' \in \cS} \rho(s') \tsum_{t=0}^\infty \gamma^t \mathtt{P}^{\pi}(s_t = s| s_0=s')$ for every $s \in \cS$, where $\mathtt{P}^{\pi}(s_t = s| s_0=s')$ denotes the probability of reaching state $s$, if running $\vartheta$ starting from $s'$ within MDP $\cM_{\PP^\pi}$.
In particular, we write $d_{s}^{\pi}$ when $\rho$ has support $\cbr{s}$. 
We next establish the difference of values for two policies of nature. 

\begin{lemma}\label{lemma_perf_diff}
For a pair of policies $\pi, \pi'$, and any $s\in \cS$,  we have
\begin{align}\label{eq_perf_diff}
\cV^{\pi'}(s) - \cV^{\pi}(s) = \frac{1}{1-\gamma}
\EE_{s' \sim d_{s}^{\pi'}} \sbr{
\cQ^{\pi}(s', \pi'(s')) - 
\cQ^{\pi}(s', \pi(s'))
}
\end{align}
Equivalently, by defining $\phi^{\pi}( s, \pi'(s)) 
%\coloneqq \cQ^{\pi}(s, \pi'(s)) - 
%\cQ^{\pi}(s, \pi(s)) 
\coloneqq \gamma \zeta \inner{\pi'(s) - \pi(s)}{\cV^{\pi}_{\vartheta, s}}$, then 
\begin{align}\label{eq_perf_diff_linearized}
\cV^{\pi'}(s) - \cV^{\pi}(s) = \frac{1}{1-\gamma}
\EE_{s' \sim d_{s}^{\pi'}} \sbr{\phi^{\pi}(s', \pi'(s'))}
\end{align}
\end{lemma}

\begin{proof}
%The proof of \eqref{eq_perf_diff} follows standard steps of performance difference lemma for finite MDPs \cite{lan2021policy, kakade2002approximately} and hence is omitted here.
%\yan{need to expand on this one}
%Let $\xi_(s)$ denote the 
Let $\xi'(s) = \cbr{s_0 = s, \pi'(s_0), s_1, \pi'(s_1), \ldots} $ denote the trajectory generated by $\pi'$ within $\mathfrak{M}$. 
That is 
\begin{align*}
s_{t+1} \sim \mathfrak{P}(\cdot|s_t, \pi'(s_t)),
\end{align*}
or equivalently, in view of \eqref{state_transit_of_nature_given_policy}, that
\begin{align}\label{state_transition_distribution_equivalence}
s_{t+1} \sim \tsum_{a \in \cA} \vartheta(a|s_t) \PP^{\pi'}_{s_t,a} (\cdot) .
\end{align}
We then obtain 
\begin{align*}
\cV^{\pi'}(s) - \cV^{\pi}(s)
& \overset{(a)}{=} \EE_{\xi'(s)} \sbr{\tsum_{t=0}^\infty \gamma^t \rbr{ \mathfrak{C}(s_t) + \gamma \cV^{\pi}(s_{t+1}) - \cV^{\pi}(s_t)}  + \cV^\pi(s_0) }   - \cV^\pi(s)  \\
& \overset{(b)}{=} \EE_{\xi'(s)} \sbr{\tsum_{t=0}^\infty \gamma^t \rbr{ \cQ^{\pi}(s_t, \pi'(s_t)) - \cV^{\pi}(s_t)}  } \\
& \overset{(c)}{=} \frac{1}{1-\gamma} \EE_{s' \sim d_s^{\pi'}} \sbr{\cQ^{\pi}(s', \pi'(s')) - \cV^\pi(s')},
\end{align*}
where $(a)$ follows from the definition of $\cV^{\pi'}(s)$, 
 $(b)$  follows from $s_0 = s$ and \eqref{relation_q_and_v},
and $(c)$ follows from \eqref{state_transition_distribution_equivalence} and the definition of $d_s^{\pi'}$.
Then \eqref{eq_perf_diff} follows  by noting that $\cV^{\pi}(s) = \cQ^\pi(s, \pi(s))$. 
Finally, \eqref{eq_perf_diff_linearized} follows from \eqref{eq_perf_diff} and Proposition \ref{prop_q_structure}.
\end{proof}

Interested readers might find the formulated MDP of nature challenging upon initial examination. 
In particular, as nature has a continuous action space, even evaluating the state-action value function \eqref{def_q_func_nature} seems to be challenging, a crucial quantity for policy improvement. 
It is also unclear whether one should and how to parameterize the policy of nature. 
In the next section, we proceed to develop the first-order robust policy evaluation (FRPE) method that exploits the structural properties established in this section and overcomes the aforementioned difficulties.

%!TEX root = ./robust_pe.tex
%\newpage

\section{Deterministic Robust Policy Evaluation}\label{sec_deterministic}
The deterministic setting we consider in this section assumes that $\overline{\PP}$ is known. 
%As a consequence,  $\PP^\pi$ defined in \eqref{kernel_defined_by_nature_policy} is available to the planner, for any fixed policy $\pi$ of nature.
We separate our discussions into two parts. 
The first part deals with the tabular setting, where the state space is relatively small, and exact computation is affordable.
We then discuss the extension to handling large state space in the second part, which involves inexact computation.

\subsection{Tabular Setting}

The proposed method, first-order robust policy evaluation (FRPE),  assumes the access to an oracle with the following capability:
for any $\pi_k$ of the nature, it returns the value function $V^{\vartheta}_{\PP^{\pi_k}}$ of $\vartheta$ within MDP $\cM_{\PP^{\pi_k}}$. 
Equipped with this oracle, FRPE then updates the policy of the nature according to 
\begin{align*}
	\pi_{k+1}(s) = \gamma \zeta \argmin_{\DD \in \cD_s}  \tsum_{t=0}^{k}  \beta_t \inner{ \DD}{{\cV}^{\pi_t}_{\vartheta, s}} + \lambda_k w(\DD), ~ \forall s \in \cS,
\end{align*}
where $w(\cdot)$ is a strictly convex function over $\Delta_{\cS}^{\abs{\cA}}$.
%where ${\cV}^{\pi_t}_{\vartheta, s} = \vartheta(\cdot|s) \otimes {\cV}^{\pi_t}$.
For tabular setting this evaluation oracle can be easily implemented. 
Given that $\overline{\PP}$ and $\pi_k$ are both known,  one can directly compute $\PP^{\pi_k}$ defined in \eqref{kernel_defined_by_nature_policy}. 
Then evaluating $V^{\vartheta}_{\PP^{\pi_k}}$ reduces to a standard policy evaluation problem with known transition kernel $\PP^{\pi_k}$, which can be computed by either solving a linear system or fixed point iteration.

\begin{algorithm}[t!]
  \caption{First-order Robust Policy Evaluation (FRPE)}
  \begin{algorithmic}
%    \REQUIRE Input
%    \ENSURE Output
    \STATE {\bf Input:} $\cbr{(\beta_k, \lambda_k)}$.
    \STATE {\bf Initialize:} arbitrary initial policy $\pi_0 \in \Pi$.
 \FOR{$k = 0, 1, \ldots$}
 	\STATE Evaluate ${\cV}^{\pi_k} = - {V}^{\vartheta}_{\PP^{\pi_k}}$.
	\STATE Update:
	\begin{align}\label{raw_update_deterministic_rpe}
	\pi_{k+1}(s) = \gamma \zeta \argmin_{\DD \in \cD_s}  \tsum_{t=0}^{k}  \beta_t \inner{ \DD}{{\cV}^{\pi_t}_{\vartheta, s}} + \lambda_k w(\DD), ~ \forall s \in \cS.
	\end{align}
	\STATE	 where 
	$
		{\cV}^{\pi_t}_{\vartheta, s} \coloneqq \vartheta(\cdot|s) \otimes {\cV}^{\pi_t} .
	$
    \ENDFOR
%    \RETURN $\cV^{\pi_k}$
     \end{algorithmic}
\end{algorithm}

It should be noted that to perform \eqref{raw_update_deterministic_rpe} one does not necessarily need to store the historical $\{\cV^{\pi_t}\}_{t=0}^k$.
Instead one can maintain a proper running average of these historical values.
From the definition of $\phi^{\pi}$ in Lemma \ref{lemma_perf_diff}, it is also clear that \eqref{raw_update_deterministic_rpe} is equivalent to the following update: 
\begin{align}
\pi_{k+1}(s) & = \argmin_{\DD \in \cD_s} \tsum_{t= 0}^k \beta_t \phi_t(s, \DD)  + \lambda_k w(\DD) \nonumber \\
& = \argmin_{\DD \in \cD_s}  \Phi_k(s, \DD)  + \lambda_k w(\DD) , \label{pda_nabular_deterministic_update}
\end{align}
where we have defined 
\begin{align*}
\phi_t \coloneqq \phi^{\pi_t}, ~ 
\Phi_k \coloneqq \tsum_{t=0}^k \beta_t \phi_t.
\end{align*} 
Going forward we will often make use of the  Bregman divergence associated with $w(\cdot)$, defined as 
\begin{align*}
\cB(\DD, \DD') \coloneqq w(\DD') - w(\DD) - \inner{\nabla  w(\DD)}{\DD' - \DD}.
\end{align*}
%We assume $w(\cdot)$ is a strictly convex function over $\Delta_{\cS}^{\abs{\cA}}$ and its associated Bregman divergence is defined as 
%\begin{align*}
%\cB(\DD, \DD') \coloneqq w(\DD) - w(\DD') - \inner{\nabla  w(\DD')}{\DD - \DD'}.
%\end{align*}
%Many of our ensuing discussions consider the case where the distance-generating function is defined as 
%\begin{align}\label{dgf_negative_entropy}
%w(\DD) =  \tsum_{a \in \cA} \tsum_{s' \in \cS} \DD_{a}(s') \log \DD_{a}(s').
%\end{align}
%In this case, it can be readily verify that  
%\begin{align}\label{bregman_divergence_negative_entropy}
%\cB(\DD, \DD') & =  \tsum_{a \in \cA} \tsum_{s' \in \cS} \DD_a(s') \log \rbr{ 
%\frac{\DD_a(s')}{\DD'_a(s)}
%} 
%% \nonumber \\ 
%% & = \tsum_{a \in \cA} \mathrm{KL}(\DD_{a} \Vert \DD'_{a})   \nonumber \\
%%& \geq \frac{1}{2\abs{\cA}} \norm{\DD - \DD'}_1^2, 
% \geq \frac{1}{2} \tsum_{a \in \cA} \norm{\DD_a - \DD'_a}_1^2,
%%~ \forall \DD, \DD' \in \Delta_{\cS}^{\abs{\cA}}.
%\end{align}
%where the last inequality follows from the Pinsker's inequality.
As \eqref{raw_update_deterministic_rpe} is invariant when shifting $w$ by a constant, without loss of generality we assume $\inf_{\DD \in \Delta_{\cS}^{\abs{\cA}}} w(\DD) \geq 0$.
We also require  $\overline{w} = \sup_{\DD \in \Delta_{\cS}^{\abs{\cA}}} w(\DD) < \infty$.

We start by  providing a simple characterization on each step of FRPE.

\begin{lemma}\label{lemma_pda_determinsitic_step_characterization}
Let $\Phi_{-1} \equiv 0$, $\lambda_{-1} = 0$, and $\lambda_k \geq \lambda_{k-1}$ for every $k \geq 1$.
Then for any $k \geq 0$, we have 
\begin{align}
& \Phi_k(s, \pi_{k+1}(s)) + \lambda_k \cB( \pi_{k+1}(s), \DD) \leq \Phi_k(s, \DD), ~ \forall \DD \in \cD_s, \label{pda_nhree_point} \\
& \beta_k \phi_k(s, \pi_{k+1}(s)) 
 \leq \Phi_k(s, \pi_{k+1}(s))  - \Phi_{k-1}(s, \pi_k(s)) - \lambda_{k-1} \cB(\pi_k(s), \pi_{k+1}(s))  \label{raw_progress_ineq_pda} .
\end{align}
%where $\cB(\DD, \DD') \coloneqq w(\DD) - w(\DD') - \inner{\nabla  w(\DD')}{\DD - \DD'}$ denotes the Bregman divergence between actions $\DD$ and $\DD'$.
\end{lemma}

\begin{proof}
Clearly, from the optimality condition of \eqref{pda_nabular_deterministic_update}, we obtain \eqref{pda_nhree_point}.
Given the definition of $\Phi_k$, we have 
\begin{align*}
\beta_k \phi_k(s, \DD) & = \Phi_k(s, \DD) - \Phi_{k-1}(s, \DD) - (\lambda_k - \lambda_{k-1}) w(\DD) \\
& \leq \Phi_k(s, \DD) - \Phi_{k-1}(s, \DD) \\
& \leq \Phi_k(s, \DD) - \Phi_{k-1}(s, \pi_k(s)) - \lambda_{k-1} \cB(\pi_k(s), \DD), ~ \forall k \geq 0.
\end{align*}
Further taking $\DD = \pi_{k+1}(s)$ in the above inequality yields \eqref{raw_progress_ineq_pda}.
\end{proof}

With Lemma \ref{lemma_pda_determinsitic_step_characterization} in place, we  proceed to establish the generic convergence properties of FRPE.

\begin{lemma}\label{lemma_pda_deterministic_generic}
Suppose $\lambda_k \geq \lambda_{k-1}$ for every $k \geq 1$, then 
\begin{align*}
\tsum_{t=0}^k 
\beta_t \rbr{f(\pi_{t+1}) - f(\pi^*)}
 \leq 
\tsum_{t=0}^k 
\beta_t \rbr{
1  - \frac{1-\gamma}{\norm{{d_{\rho}^{\pi^*}}/{\rho}}_\infty}
}
\rbr{f(\pi_t) - f(\pi^*)} 
+  \frac{2 \lambda_k\overline{w}}{1-\gamma}.
\end{align*}
%where $\overline{w} = \sup_{\DD \in \Delta_{\cS}^{\abs{\cA}}} w(\DD) $. 
\end{lemma}

\begin{proof}
%Clearly, from the optimality condition of \eqref{pda_nabular_deterministic_update}, we obtain
%\begin{align}
%\label{pda_nhree_point}
%\Phi_k(s, \pi_{k+1}(s)) + \lambda_k \cB(\DD, \pi_{k+1}(s)) \leq \Phi_k(s, \DD), ~ \forall \DD \in \cD_s,
%\end{align} 
%where $\cB(\DD, \DD') = w(\DD) - w(\DD') - \inner{\nabla  w(\DD')}{\DD - \DD'}$ denotes the Bregman divergence between actions $\DD$ and $\DD'$.
%Given the definition of $\Phi_k$, we have 
%\begin{align*}
%\beta_k \phi_k(s, \DD) & = \Phi_k(s, \DD) - \Phi_{k-1}(s, \DD) - (\lambda_k - \lambda_{k-1}) w(\DD) \\
%& \leq \Phi_k(s, \DD) - \Phi_{k-1}(s, \DD) \\
%& \leq \Phi_k(s, \DD) - \Phi_{k-1}(s, \pi_k(s)) - \lambda_{k-1} \cB(\DD, \pi_k(s)), ~ \forall k \geq 0,
%\end{align*}
%where we have defined $\Phi_{-1} \equiv 0$ and $\lambda_{-1} = 0$.
%Further taking $\DD = \pi_{k+1}(s)$ in the above inequality yields 
From \eqref{raw_progress_ineq_pda} in Lemma \ref{lemma_pda_determinsitic_step_characterization} we have
\begin{align}
\beta_k \phi_k(s, \pi_{k+1}(s)) 
& \leq \Phi_k(s, \pi_{k+1}(s))  - \Phi_{k-1}(s, \pi_k(s)) - \lambda_{k-1} \cB( \pi_k(s), \pi_{k+1}(s)) \nonumber \\
& \overset{(a)}{\leq} \Phi_k(s, \pi_{k}(s)) - \Phi_{k-1}(s, \pi_k(s)) - \lambda_{k-1} \cB(\pi_k(s), \pi_{k+1}(s))  - \lambda_k \cB( \pi_{k+1}(s), \pi_k(s))  \nonumber  \\
& \overset{(b)}\leq \beta_k \phi_k(s, \pi_k(s)) + (\lambda_k - \lambda_{k-1}) w(\pi_k(s)) \nonumber \\
&\overset{(c)} = (\lambda_k - \lambda_{k-1}) w(\pi_k(s)), \nonumber \\
& \leq (\lambda_k - \lambda_{k-1}) \overline{w} \label{subgrad_inner_ub}, 
\end{align}
where $(a)$ follows from applying \eqref{pda_nhree_point} again;
$(b)$ follows from the definition of $\Phi_k$;
 $(c)$ follows from the trivial identity $\phi_k(s, \pi_k(s)) = 0$ given the definition of $\phi_k$.
 Combing \eqref{subgrad_inner_ub} with Lemma \ref{lemma_perf_diff}, we have
\begin{align}
\cV^{\pi_{k+1}}(s) - \cV^{\pi_k}(s) 
& = \frac{1}{1-\gamma} \EE_{d_s^{\pi_{k+1}}} \sbr{
\phi_k(s', \pi_{k+1}(s') )- \frac{(\lambda_k - \lambda_{k-1}) \overline{w}}{\beta_k} 
} 
+  \frac{(\lambda_k - \lambda_{k-1} ) \overline{w}}{ (1-\gamma) \beta_k } 
 \nonumber \\
& \leq 
\frac{1}{1-\gamma} d_s^{\pi_{k+1}}(s) \sbr{
\phi_k(s', \pi_{k+1}(s') )- \frac{(\lambda_k - \lambda_{k-1}) \overline{w}}{\beta_k} 
}
+ \frac{(\lambda_k - \lambda_{k-1} ) \overline{w}}{ (1-\gamma) \beta_k } \nonumber  \\
& = 
\phi_k(s, \pi_{k+1}(s) ) +  \frac{ \gamma(\lambda_k - \lambda_{k-1}) \overline{w}}{(1-\gamma) \beta_k} . \label{eq_pda_approx_progress}
\end{align}
It is also clear from the first equality above and \eqref{subgrad_inner_ub} that 
\begin{align}
\cV^{\pi_{k+1}}(s) - \cV^{\pi_k}(s) \leq  \frac{ (\lambda_k - \lambda_{k-1}) \overline{w}}{(1-\gamma) \beta_k} . \label{eq_pda_approx_monotone}
\end{align}
Now taking the telescopic sum of \eqref{raw_progress_ineq_pda}, we obtain 
\begin{align*}
\tsum_{t=0}^k \beta_t \phi_t(s, \pi_{t+1}(s))
& \leq \Phi_k(s, \pi_{k+1}(s)) - \Phi_{-1}(s, \pi_0(s)) \\
%- \lambda_{-1} \cB(\pi_0(s), \pi_1(s)) \\
& \overset{(d)}{=} \Phi_k(s, \pi_{k+1}(s)) \\
& \overset{(e)}{\leq} \Phi_k(s, \DD) = \tsum_{t=0}^k \beta_t \phi_t(s, \DD) + \lambda_k w(\DD) .
\end{align*}
where $(d)$ applies the definition of $\Phi_{-1} \equiv 0$,
and $(e)$ applies \eqref{pda_nhree_point}.
Combining the above inequality with \eqref{eq_pda_approx_progress}, 
we have 
\begin{align*}
\tsum_{t=0}^k \beta_t \rbr{ \cV^{\pi_{t+1}}(s) - \cV^{\pi_t}(s) -  \frac{(\lambda_t - \lambda_{t-1} ) \overline{w}}{\beta_t (1-\gamma)}}
\leq \tsum_{t=0}^k \beta_t \phi_t(s, \DD)  +  \lambda_k \overline{w}.
\end{align*}
Substituting  $\DD = \pi^*(s)$ into the above for every $s\in \cS$,  and aggregating the resulting inequalities by weight $d_{s }^{\pi^*}$, we obtain 
\begin{align*}
& \tsum_{t=0}^k \EE_{s' \sim d_{s}^{\pi^*}} \big[   \beta_t ( \cV^{\pi_{t+1}}(s) - \cV^{\pi_t}(s) -  \frac{(\lambda_t - \lambda_{t-1} ) \overline{w}}{\beta_t (1-\gamma)}) \big] \\
 \leq &  \tsum_{t=0}^k \beta_t \EE_{s' \sim d_s^{\pi^*}}\sbr{ \phi_t(s, \pi^*(s))} +  \lambda_k \overline{w} \\
 \overset{(f)}{=} & 
(1-\gamma) \tsum_{t=0}^k \beta_t \rbr{\cV^{\pi^*}(s) - \cV^{\pi_t}(s)} + \lambda_k \overline{w},
\end{align*} 
where $(f)$ applies Lemma \ref{lemma_perf_diff}.
It remains to take expectation with respect to $s \sim \rho$ of the above inequality and obtain
%Combining the above inequality and  \eqref{eq_pda_approx_monotone}, we can further obtain 
\begin{align*}
&  \tsum_{t=0}^k  \big\lVert\frac{d_{\rho}^{\pi^*}}{\rho}\big\rVert_\infty \EE_{s \sim \rho}
\sbr{ \beta_t ( \cV^{\pi_{t+1}}(s) - \cV^{\pi_t}(s) -  \frac{(\lambda_t - \lambda_{t-1} ) \overline{w}}{\beta_t (1-\gamma)})} \\
\overset{(h)}{\leq} & 
\tsum_{t=0}^k \EE_{s \sim d_{\rho}^{\pi^*}} \sbr{  \beta_t ( \cV^{\pi_{t+1}}(s) - \cV^{\pi_t}(s) -  \frac{(\lambda_t - \lambda_{t-1} ) \overline{w}}{\beta_t (1-\gamma)})} \\
 \leq & 
 (1-\gamma) \tsum_{t=0}^k \beta_t \EE_{s \sim \rho} \sbr{\cV^{\pi^*}(s) - \cV^{\pi_t}(s)} + \lambda_k \overline{w},
\end{align*} 
where $(h)$ utilizes \eqref{eq_pda_approx_monotone}, from which the term inside expectation is non-positive.
Simple arrangement of the above relation yields 
\begin{align*}
\tsum_{t=0}^k 
\beta_t \rbr{f(\pi_{t+1}) - f(\pi^*)}
& \leq 
\tsum_{t=0}^k 
\beta_t \rbr{
1  - \frac{1-\gamma}{\norm{{d_{\rho}^{\pi^*}}/{\rho}}_\infty}
}
\rbr{f(\pi_t) - f(\pi^*)} 
+  \rbr{ \frac{1}{1-\gamma} + {\norm{{d_{\rho}^{\pi^*}}/{\rho}}_\infty}^{-1} } \lambda_k\overline{w}  \\
& \leq 
\tsum_{t=0}^k 
\beta_t \rbr{
1  - \frac{1-\gamma}{\norm{{d_{\rho}^{\pi^*}}/{\rho}}_\infty}
}
\rbr{f(\pi_t) - f(\pi^*)} 
+ \frac{ 2 \lambda_k \overline{w}}{1-\gamma}  .
\end{align*}
The proof is then completed.
\end{proof}

We are ready to establish the linear convergence of FRPE with proper specification of $\cbr{(\beta_k, \lambda_k)}$.

\begin{theorem}\label{thrm_linear_convergence_frpe}
Take $\beta_k = \rbr{
1  - \frac{1-\gamma}{\norm{{d_{\rho}^{\pi^*}}/{\rho}}_\infty}
}^{-k}$ and $\lambda_k = \lambda \zeta > 0 $,  then 
\begin{align*}
f(\pi_{k+1}) - f(\pi^*) \leq 
\rbr{
1  - \frac{1-\gamma}{\norm{{d_{\rho}^{\pi^*}}/{\rho}}_\infty}
}^{k} \sbr{ f(\pi_0) - f(\pi^*) + \frac{2 \lambda \zeta \overline{w}}{1-\gamma} }.
\end{align*}
\end{theorem}

\begin{proof}
The  claim follows from a direct application of Lemma \ref{lemma_pda_deterministic_generic}, combined with the choice of $\cbr{(\beta_k, \lambda_k)}$.
\end{proof}

%\yan{remark on comparison to standard approach using fixed point iteration, mention the failure of the latter approach to function approximation, and then motivate the discussion of linear function approximation}

In view of Theorem \ref{thrm_linear_convergence_frpe}, FRPE in the deterministic setting attains linear convergence.
In particular, the performance of applying FRPE to the offline robust MDP problem described in Example \ref{example_offline_rmdp}  matches the classical robust policy evaluation methods based on fixed point iteration \cite{nilim2005robust, iyengar2005robust}. 
Notably both methods requires computing the robust value for every state. 
In the following section we proceed to discuss the setting where approximate computation is required, and demonstrate the clear advantage of FRPE.

%Now suppose 
%\begin{align*}
%\beta_t \rbr{
%1  - \frac{1-\gamma}{\norm{{d_{\rho}^{\pi^*}}/{\rho}}_\infty}
%}
%\leq \beta_{n-1}, ~ \beta_0 = 1,
%\end{align*}
%we then obtain 
%\begin{align*}
%\beta_k \rbr{f(\pi_{k+1}) - f(\pi^*)}
%\leq   f(\pi_0) - f(\pi^*) + 2 \lambda_k \overline{w},
%\end{align*}
%which implies 
%\begin{align*}
%f(\pi_{k+1}) - f(\pi^*) \leq 
%\rbr{
%1  - \frac{1-\gamma}{\norm{{d_{\rho}^{\pi^*}}/{\rho}}_\infty}
%}^{k} \sbr{ f(\pi_0) - f(\pi^*) + 2 \lambda \overline{w} }
%\end{align*}

%
%\yan{add the constant stepsize version, show its iteration complexity scales with $\zeta$, hence more favorable when the ambiguity set's size $\zeta$ becomes small
%
%Wait, it seems that the iteration complexity above already implicitly scales with $\zeta$, to see this, simply take $\lambda = \zeta$, and the right hand side both scales as $\cO(\zeta)$.
%}

\subsection{FRPE with Function Approximation}

We now discuss the extension of FRPE to handle large state space.
In such a setting both parameterizing the policy $\pi$ of the nature and its value function $\cV^{\pi}$ can be costly.
We start by observing the that FRPE does not require explicit parameterization  of $\pi$.
Namely, the latest policy $\pi_{t+1}(s)$ can be generated incrementally according to \eqref{raw_update_deterministic_rpe} whenever its value is needed. 
On the other hand, parameterizing $\cV^{\pi}$ can be efficiently done by employing function approximation with relatively few parameters. 
We now discuss an example of such approach, and establish the convergence of FRPE in the presence of potential approximation error.

\begin{example}
Suppose a feature map $\psi: \cS \to \RR^d$ is given, we consider the problem of approximating $\cV^{\pi}(\cdot)$ by $\psi(\cdot)^\top \theta^\pi$ so the difference of these two functions can be properly controlled. 
In view of Lemma \ref{lemma_value_correspondence}, this is equivalent to approximating $V^{\vartheta}_{\PP^\pi}(s)$ by $-\psi(\cdot)^\top \theta^\pi$, a long studied problem with numerous off-the-shelf methods.
%which by itself is simply approximating the value function of $\vartheta$ with respect to MDP $\cM_{\PP^\pi}$.
%Such a problem has been long studied in the reinforcement learning and dynamic programming literature. 
In particular, 
as $\PP^\pi$ is known,
one can consider solving the so-called projected Bellman equation
\begin{align}\label{projected_bellman_equation}
\Psi \theta^\pi = \Pi_{\Psi, \nu^{\vartheta}_{\PP^\pi}} \cT (\Psi \theta^\pi),
\end{align}
where $\Psi \in \RR^{\abs{\cS} \times d}$ denotes the feature matrix constructed by applying $\psi(\cdot)$ to every state,
$\Pi_{\Psi, \nu^{\vartheta}_{\PP^\pi}}$ denotes the projection onto $\mathrm{span}(\Psi)$ in $\norm{\cdot}_{\nu^{\vartheta}_{\PP^\pi}}$ norm,
$ \nu^{\vartheta}_{\PP^\pi}$ denotes the stationary state distribution induced by $\vartheta$ within MDP $\cM_{\PP^\pi}$,
and $\cT$ denotes the Bellman operator of $\vartheta$ within MDP $\cM_{\PP^\pi}$.
A feature to note here is that solution $\theta^\pi$ exists for \eqref{projected_bellman_equation} \cite{580874}, and the latter  being a linear system of $\theta^\pi$ implies efficient methods for solving it \cite{karczmarz1937angenaherte, li2023accelerated}.
The approximation error is captured by $\norm{\Psi \theta^\pi - \cV^\pi}_\infty$, which depends on the expressiveness of $\psi$. 
An alternative to solving \eqref{projected_bellman_equation} is to minimize the mean-squared Bellman error, which also yields a linear system and does not require information of $\nu^{\vartheta}_{\PP^{\pi}}$.
In Section \ref{sec_stoch_linear_approx} we study this approach in the stochastic~setting.
\end{example}

In the remainder of this section we proceed to establish the convergence of FRPE with arbitrary approximation scheme, as long as the error can be controlled. 
%We now turn our attention to the convergence of FRPE where potential function approximation error exists in the evaluated $\cV^\pi$. 
Namely, instead of using exact $\cV^{\pi_t}$ in the update \eqref{raw_update_deterministic_rpe} of FRPE, one only has an approximation  $\tilde{\cV}^{\pi_t}$ such that 
\begin{align}\label{funx_approx_deterministic}
\norm{\cV^{\pi_t} - \tilde{\cV}^{\pi_t} }_\infty \leq \varepsilon_{\mathrm{approx}} 
\end{align}
for some $\varepsilon_{\mathrm{approx}}  > 0$.
%This substantially generalizes prior approaches mentioned above, by lifting any technical conditions and being applicable to broader approximation schemes.
The FRPE method in this setting takes the following update:
\begin{align}\label{raw_update_deterministic_rpe_func_approx}
\textstyle
	\pi_{k+1}(s) = \gamma \zeta \argmin_{\DD \in \cD_s}  \tsum_{t=0}^{k}  \beta_t \inner{ \DD}{\tilde{\cV}^{\pi_t}_{\vartheta, s}} + \lambda_k w(\DD), ~ \forall s \in \cS,
\end{align}
where $\tilde{\cV}^{\pi_t}_{\vartheta, s} \coloneqq \vartheta(\cdot|s) \otimes \tilde{\cV}^{\pi_t}$.
Clearly, \eqref{raw_update_deterministic_rpe_func_approx} is similar to \eqref{raw_update_deterministic_rpe} except that we replace the exact value function $\tilde{\cV}^{\pi_t}$ with its approximation $\tilde{\cV}^{\pi_t}$.
%We proceed to show that FRPE remains stable in the presence of approximation error.
It is clear that  \eqref{raw_update_deterministic_rpe_func_approx} is equivalent to:
\begin{align}
\pi_{k+1}(s) & = \argmin_{\DD \in \cD_s} \tsum_{t= 0}^k \beta_t \tilde{\phi}_t(s, \DD)  + \lambda_k w_s(\DD) \nonumber \\
& = \argmin_{\DD \in \cD_s}  \tilde{\Phi}_k(s, \DD)  + \lambda_k w_s(\DD) , \label{pda_nabular_func_approx_update}
\end{align}
where  
$
\tilde{\phi}_t(s, \DD)  = \gamma \zeta \inner{\DD - \pi_t(s)}{\tilde{\cV}^{\pi_t}_{\vartheta, s}}$ and  $\tilde{\Phi}_k  = \tsum_{t=0}^k \beta_t \tilde{\phi}_t .
$
%Let us also define 
%\begin{align}\label{def_noise_in_phi}
%\delta_n(s, \DD)  \coloneqq \tilde{\phi}_t(s, \DD)  - \phi_t (s, \DD) = \gamma \zeta \inner{\DD - \pi_t(s)}{\tilde{\cV}^{\pi}_{\vartheta, s} - {\cV}^{\pi}_{\vartheta, s}},
%\end{align}
%where the last equality follows from the definition of $\tilde{\phi}_t$ and $\phi_t$.

The following lemma follows from similar lines as in Lemma \ref{lemma_pda_determinsitic_step_characterization}.

\begin{lemma}\label{lemma_pda_determinsitic_step_characterization_func_approx}
Let $\Phi_{-1} \equiv 0$, $\lambda_{-1} = 0$, and $\lambda_k \geq \lambda_{k-1}$ for every $k \geq 1$.
Then for any $k \geq 0$, we have 
\begin{align}
& \tilde{\Phi}_k(s, \pi_{k+1}(s)) + \lambda_k \cB(\pi_{k+1}(s), \DD) \leq \tilde{\Phi}_k(s, \DD), ~ \forall \DD \in \cD_s, \label{pda_nhree_point_func_approx} \\
& \beta_k \tilde{\phi}_k(s, \pi_{k+1}(s)) 
 \leq \tilde{\Phi}_k(s, \pi_{k+1}(s))  - \tilde{\Phi}_{k-1}(s, \pi_k(s)) - \lambda_{k-1} \cB(\pi_k(s), \pi_{k+1}(s))  \label{raw_progress_ineq_pda_func_approx} .
\end{align}
%where $\cB(\DD, \DD') \coloneqq w(\DD) - w(\DD') - \inner{\nabla  w(\DD')}{\DD - \DD'}$ denotes the Bregman divergence between actions $\DD$ and $\DD'$.
\end{lemma}

We now proceed to establish the convergence properties of FRPE in the presence of approximation error.

\begin{theorem}
Suppose \eqref{funx_approx_deterministic} holds, 
then taking $\beta_k = \rbr{
1  - \frac{1-\gamma}{\norm{{d_{\rho}^{\pi^*}}/{\rho}}_\infty}
}^{-k}$ and $\lambda_k = \lambda \zeta > 0 $ yields 
\begin{align*}
f(\pi_{k+1}) - f(\pi^*) \leq 
\rbr{
1  - \frac{1-\gamma}{\norm{{d_{\rho}^{\pi^*}}/{\rho}}_\infty}
}^{k} \sbr{ f(\pi_0) - f(\pi^*) + \frac{2 \lambda \zeta \overline{w}}{1-\gamma} }
+ \frac{2 (1+ \norm{{d_\rho^{\pi^*}}/{\rho}}_\infty) \gamma \zeta  \varepsilon_{\mathrm{approx}} }{(1-\gamma)^2}.
\end{align*}
\end{theorem}

\begin{proof}
We start by noting that following the same lines as we show \eqref{subgrad_inner_ub}, one has 
\begin{align}
\beta_k \tilde{\phi}_k(s, \pi_{k+1}(s)) 
 \leq (\lambda_k - \lambda_{k-1}) \overline{w} \label{subgrad_inner_ub_func_approx}.
\end{align}
Similar to \eqref{eq_pda_approx_progress}, we obtain 
\begin{align}
& \cV^{\pi_{k+1}} (s) - \cV^{\pi_k} (s) 
-  \frac{1}{1-\gamma} \EE_{s' \sim d_s^{\pi_{k+1}}} \sbr{
\phi_k(s', \pi_{k+1}(s')) - \tilde{\phi}_k(s', \pi_{k+1}(s')) 
} \nonumber \\
\overset{(a)}{ =} & \frac{1}{1-\gamma} \EE_{s' \sim d_s^{\pi_{k+1}}} \sbr{
\tilde{\phi}_k(s', \pi_{k+1}(s')) 
}  \nonumber \\
= & \frac{1}{1-\gamma} \EE_{s' \sim d_s^{\pi_{k+1}}} \sbr{
\tilde{\phi}_k(s', \pi_{k+1}(s')) - \frac{(\lambda_k - \lambda_{k-1})\overline{w}}{\beta_k}
}
+ \frac{(\lambda_k - \lambda_{k-1}) \overline{w}}{(1-\gamma)\beta_k}  \label{approx_progress_func_approx_raw_1} \\
\overset{(b)}{\leq} & 
\tilde{\phi}_k(s, \pi_{k+1}(s)) + \frac{(\lambda_k - \lambda_{k-1} )\overline{w}}{\beta_k (1-\gamma)}, \label{inner_product_lb_func_approx}
\end{align}
where $(a)$ follows from Lemma \ref{lemma_perf_diff}, and $(b)$ follows from \eqref{subgrad_inner_ub_func_approx}.
In addition, from the definition of $\phi_k$ and $\tilde{\phi}_k$, one has 
\begin{align}\label{func_approx_error_effect_on_q}
\abs{\phi_k(s, \DD) - \tilde{\phi}_k(s,\DD)}
& = \abs{
\tsum_{a \in \cA} \gamma \zeta \inner{\DD_a - \DD^{\pi_k(s)}_a}{ \cV^{\pi_k} - \hat{\cV}^{\pi_k} } \vartheta(a|s)
}  \nonumber \\
& \leq 2 \gamma \zeta \norm{\cV^{\pi_k} - \hat{\cV}^{\pi_k}}_\infty \leq 2  \gamma \zeta \varepsilon_{\mathrm{approx}}.
\end{align}
where the last inequality follows from H\"{o}lder's inequality.
Combining \eqref{subgrad_inner_ub_func_approx},  \eqref{approx_progress_func_approx_raw_1},  and \eqref{func_approx_error_effect_on_q} also yields 
\begin{align}\label{approx_progress_func_approx_raw_3}
 \cV^{\pi_{k+1}} (s) - \cV^{\pi_k} (s)  \leq \frac{(\lambda_k - \lambda_{k-1}) \overline{w}}{(1-\gamma) \beta_k} + \frac{2 \gamma \zeta \varepsilon_{\mathrm{approx}}}{1-\gamma}.
\end{align}
Repeating the same lines after \eqref{eq_pda_approx_monotone} in the proof of Lemma \ref{lemma_pda_deterministic_generic},
with \eqref{raw_progress_ineq_pda_func_approx}, and  \eqref{inner_product_lb_func_approx}, \eqref{approx_progress_func_approx_raw_3} replaced by \eqref{raw_progress_ineq_pda}, 
\eqref{eq_pda_approx_progress}, and \eqref{eq_pda_approx_monotone},
we obtain 
\begin{align*}
\tsum_{t=0}^k 
\beta_t \rbr{f(\pi_{t+1}) - f(\pi^*)}
& \leq 
\tsum_{t=0}^k 
\beta_t \rbr{
1  - \frac{1-\gamma}{\norm{{d_{\rho}^{\pi^*}}/{\rho}}_\infty}
}
\rbr{f(\pi_t) - f(\pi^*)} 
+  \frac{2 \lambda_k\overline{w}}{1-\gamma} \\
& ~~~~~~
+ \tsum_{t=0}^k \beta_t \frac{2 \gamma \zeta  \varepsilon_{\mathrm{approx}}}{1-\gamma} \rbr{1 + \norm{\frac{d_\rho^{\pi^*}}{\rho}}_\infty^{-1}}.
\end{align*}
Plugging the choice of $\cbr{(\beta_k, \lambda_k)}$ yields the desired result.
%\begin{align*}
%f(\pi_{k+1}) - f(\pi^*) \leq 
%\rbr{
%1  - \frac{1-\gamma}{\norm{{d_{\rho}^{\pi^*}}/{\rho}}_\infty}
%}^{k} \sbr{ f(\pi_0) - f(\pi^*) + \frac{2 \lambda \zeta \overline{w}}{1-\gamma} }
%+ \frac{2 (1+ \norm{{d_\rho^{\pi^*}}/{\rho}}_\infty) \gamma \zeta  \varepsilon_{\mathrm{approx}} }{(1-\gamma)^2}.
%\end{align*}
\end{proof}

%\yan{comparison to Mannor's prior work after technical discussion}
%An immediate observation is that FRPE does not necessarily require parameterizing the policy. 
%Instead, the updated policy $\pi_{t+1}(s)$ at any given state $s$ can be computed whenever it is need within the evaluation procedure. 

At this point it is worth mentioning a few prior studies that take a direct approach towards robust policy evaluation with linear function approximation, by extending the projected Bellman equation \eqref{projected_bellman_equation} to the robust setting \cite{roy2017reinforcement, tamar2014scaling}.
In this case, the objective is to solve 
\begin{align}\label{robust_projected_bellman_equation}
\Psi \theta^* = \Pi_{\Psi,  \nu} \cT_{\mathrm{r}} (\Psi \theta^*)
\end{align}
 to obtain estimate $\Psi \theta^* \approx V^{\vartheta}_r$. 
Compared to \eqref{projected_bellman_equation}, the operator $\cT_{\mathrm{r}}$ corresponds to the so-called robust Bellman operator defined in \eqref{def_robust_ballmen_op}, and $\nu$ is the stationary distribution for certain exploration policy.
An important limitation of such approach is that \eqref{robust_projected_bellman_equation} does not necessarily admit a solution, as the operator $\Pi_{\Psi,  \nu} \cT_{\mathrm{r}}$ is no longer a contraction. 
Indeed it is shown in \cite{tamar2014scaling} that an restrictive yet necessary condition is required to certify the existence of $\theta^*$ satisfying \eqref{robust_projected_bellman_equation}.
This assumption also seems difficult to verify even if the model is known to the planner, and only asymptotic convergence has been established depsite its restrictive nature.

On the other hand, if one seeks the approximate solution of \eqref{robust_projected_bellman_equation} by minimizing its mean-squared error $\norm{\Psi \theta^* - \Pi_{\Psi,  \nu} \cT_{\mathrm{r}} (\Psi \theta^*)}^2_2$, the resulting objective can be easily non-convex due to the non-linearity of $\cT_{\mathrm{r}}$. 
 It should be clear that FRPE substantially improves the aforementioned approaches, by removing unrealistic assumptions while being applicable to broader approximation schemes.

It should also be noted that for large-scale offline robust MDPs, both the model $\overline{\PP}$ and ambiguity set $\cD$ can be difficult to store. 
Consequently, in solving the linear system defined by \eqref{projected_bellman_equation} or the mean-squared Bellman error, one might proceed with an incremental manner.   
This can be done, for instance, by Kaczmarz method \cite{karczmarz1937angenaherte} and its randomized variants \cite{strohmer2009randomized, gower2015randomized}.
One can also utilize the stochastic variant of FRPE to be discussed in the next section.

%\yan{should mention that only **asymptotic results** are obtained with this approach!!!}
%
%\yan{should we mention that stochasticity is also allowed even when $\overline{\PP}$ is known }

%!TEX root = ./robust_pe.tex

\section{Stochastic Robust Policy Evaluation}\label{sec_stochastic}

Compared to the deterministic setting, in the stochastic setting we do not have exact information on $\overline{\PP}$. Instead only transition samples with distribution governed by $\overline{\PP}$ are available.
As in Section \ref{sec_deterministic}, we initiate our discussion with tabular setting.
%, where the size of the state space is relatively small compared to available computational resource. 
Then we proceed to discuss the convergence of the proposed method with linear function approximation to handle large state spaces.\footnote{
The method to be developed can indeed be combined with general function approximation schemes, provided certain off-policy evaluation problem can be solved to a properly pre-specified accuracy. We defer related discussions to Section~\ref{sec_conclusion}. 
}
In both cases we present sample complexities that appears to be completely new for stochastic robust policy evaluation.

We start by introducing the framework of stochastic first-order robust policy evaluation (SFRPE).
Compared to Section \ref{sec_deterministic}, 
the distance generating function $w_s(\cdot)$ in the stochastic setting can potentially depend on $s$ through $\vartheta(\cdot|s)$. 
We assume in addition that the Bregman divergence induced by $w_s(\cdot)$ satisfies 
\begin{align}\label{sc_in_group_norm}
\cB_s(\DD, \DD') \coloneqq w_s(\DD') - w_s(\DD) - \inner{\nabla  w_s(\DD)}{\DD' - \DD}
 \geq  \frac{\mu_w}{2} \tsum_{a \in \cA} \vartheta(a|s) \norm{\DD_a - \DD'_a}_1^2
\end{align}
for some $\mu_w > 0$.
Many of our ensuing discussions consider the case where   
\begin{align}\label{dgf_negative_entropy}
w_s(\DD) =  \tsum_{a \in \cA} \vartheta(a|s) \tsum_{s' \in \cS} \DD_{a}(s') \log \DD_{a}(s') + \log \abs{\cS}.
\end{align}
In this case, it can be readily verified that $\cB_s(\DD, \DD')  =  \tsum_{a \in \cA}\vartheta(a|s) \tsum_{s' \in \cS}  \DD'_a(s') \log \rbr{ 
{\DD'_a(s')}/{\DD_a(s')}
}$,
and 
\begin{align}
0  \leq   w_s(\DD)  & \leq \log \abs{\cS}, \label{bregman_divergence_negative_entropy_lb_ub} \\
\cB_s(\DD, \DD')  
%=  \tsum_{a \in \cA}\vartheta(a|s) \tsum_{s' \in \cS}  \DD_a(s') \log \rbr{ 
%\frac{\DD_a(s')}{\DD'_a(s')}
%} 
&  \geq \frac{1}{2} \tsum_{a \in \cA}  \vartheta(a|s) \norm{\DD_a - \DD'_a}_1^2, \label{bregman_divergence_negative_entropy_sc}
%~ \forall \DD, \DD' \in \Delta_{\cS}^{\abs{\cA}}.
\end{align}
where \eqref{bregman_divergence_negative_entropy_sc} follows from the Pinsker's inequality.
Similar to the deterministic setting, we require $\overline{w} = \sup_{s \in \cS} \sup_{\DD \in \Delta_{\cS}^{\abs{\cA}}} w_s(\DD) < \infty$.
In view of \eqref{bregman_divergence_negative_entropy_lb_ub} and \eqref{bregman_divergence_negative_entropy_sc}, for $w_s(\cdot)$ defined in \eqref{dgf_negative_entropy} one can take $\overline{w} =\log \abs{\cA}$ and $\mu_w = 1$.

\begin{algorithm}[t]
  \caption{Stochastic First-order Robust Policy Evaluation (SFRPE)}
  \begin{algorithmic}
%    \REQUIRE Input
%    \ENSURE Output
    \STATE {\bf Input:} $\cbr{(\beta_k, \lambda_k)}$.
    \STATE {\bf Initialize:} arbitrary initial policy $\pi_0 \in \Pi$.
    \FOR{$k = 0, 1, \ldots$}
 	\STATE  Run stochastic evaluation operator to obtain $\hat{\cV}^{\pi_k} $.
%	form $\hat{\cV}^{\pi_n}_{\vartheta, s}$ for every $s \in \cS$, where 
%	\begin{align*}
%	\hat{\cV}^{\pi_n}_{\vartheta, s} \coloneqq \vartheta(\cdot|s) \otimes \hat{\cV}^{\pi_n} .
%	\end{align*}
	\STATE  Update:
	\begin{align}\label{raw_update_stoch_rpe}
	\textstyle
	\pi_{k+1}(s) = \gamma \zeta \argmin_{\DD \in \cD_s}  \tsum_{t=0}^{k}  \beta_t \inner{ \DD}{\hat{\cV}^{\pi_t}_{\vartheta, s}} + \lambda_k w_s(\DD), ~ \forall s \in \cS,
	\end{align}
	where 
	$
		\hat{\cV}^{\pi_t}_{\vartheta, s} \coloneqq \vartheta(\cdot|s) \otimes \hat{\cV}^{\pi_t} .
	$
    \ENDFOR
    \RETURN $\tsum_{t=1}^k \theta_t \hat{\cV}^{\pi_t}$, where 
    \begin{align}\label{def_theta}
	\theta_t = \beta_t / (\tsum_{t=1}^k \beta_t) .
	\end{align}
     \end{algorithmic}
\end{algorithm}

\begin{remark}
%It turns out that divergence \eqref{dgf_negative_entropy} enjoys at least two apparent advantages. 
%In purely online robust MDP problems (cf. Example \ref{online_robust_mdp}) where $\cD_s = \Delta_{\cS}^{\abs{\cA}}$, it can readily seen that update of SFRPE \eqref{raw_update_stoch_rpe} has closed-form solution
%$
%\DD^{\pi_{n+1}(s)}_a \propto \exp (
%- \gamma \zeta \tsum_{t=0}^n \hat{\cV}^{\pi_t} / \lambda_n
%)
%$
%for any $a \in \cA$.
The weighted construction of divergence \eqref{dgf_negative_entropy}  appears to be essential. In particular, setting equal weights in \eqref{dgf_negative_entropy} would result in a sample complexity that linearly depends on the action space $\abs{\cA}$ despite we are evaluating the robust value function.
%
%\yan{Need to mention on the motivation of state-dependent divergence -- remove the dependence on the action space}
%\yan{need to mention how the update can be solved efficiently for large state space, or has explicit solution}
%\yan{both are related to the choice of divergence}
\end{remark}

It is clear that update \eqref{raw_update_stoch_rpe} is equivalent to the following update
\begin{align}
\pi_{k+1}(s) & = \argmin_{\DD \in \cD_s} \tsum_{t= 0}^k \beta_t \hat{\phi}_t(s, \DD)  + \lambda_k w_s(\DD) \nonumber \\
& = \argmin_{\DD \in \cD_s}  \hat{\Phi}_k(s, \DD)  + \lambda_k w_s(\DD) , \label{pda_tabular_stoch_update}
\end{align}
where  
$\hat{\phi}_t(s, \DD)  = \gamma \zeta \inner{\DD - \pi_t(s)}{\hat{\cV}^{\pi_t}_{\vartheta, s}}$ and $ \hat{\Phi}_k  = \tsum_{t=0}^k \beta_t \hat{\phi}_t$.
Let us also define 
\begin{align}\label{def_noise_in_phi}
\delta_t(s, \DD)  \coloneqq \hat{\phi}_t(s, \DD)  - \phi_t (s, \DD) = \gamma \zeta \inner{\DD - \pi_t(s)}{\hat{\cV}^{\pi_t}_{\vartheta, s} - {\cV}^{\pi_t}_{\vartheta, s}},
\end{align}
where the last equality follows from the definition of $\hat{\phi}_t$ and $\phi_t$.

The following lemma follows the exact same argument as in Lemma \ref{lemma_pda_determinsitic_step_characterization}.

\begin{lemma}\label{lemma_pda_stochastic_step_characterization}
Define $\hat{\Phi}_{-1} \equiv 0$ and $\lambda_{-1} = 0$, and let $\lambda_k \geq \lambda_{k-1}$ for every $k \geq 1$.
Then for any $k \geq 0$, we have 
\begin{align}
& \hat{\Phi}_k(s, \pi_{k+1}(s)) + \lambda_k \cB_s(\pi_{k+1}(s), \DD) \leq \hat{\Phi}_k(s, \DD), ~ \forall \DD \in \cD_s, \label{pda_nhree_point_stoch} \\
& \beta_k \hat{\phi}_k(s, \pi_{k+1}(s)) 
 \leq \hat{\Phi}_k(s, \pi_{k+1}(s))  - \hat{\Phi}_{k-1}(s, \pi_k(s)) - \lambda_{k-1} \cB_s( \pi_k(s), \pi_{k+1}(s))  \label{raw_progress_ineq_pda_stoch} .
\end{align}
\end{lemma}

We next establish some generic convergence properties of SFRPE.

\begin{lemma}\label{lemma_generic_prop_stoch}
Let  $\lambda_t \geq \lambda_{t-1}$ for every $n \geq 1$,  $\beta_0 = 0$.
%For any $k \geq 1$, define 
%\begin{align}\label{def_theta}
%\theta_t = \beta_t / (\tsum_{t=1}^k \beta_t) .
%\end{align}
Then
%Suppose $\EE_{|k} \sbr{\norm{\hat{\cV}^{\pi_k}}_\infty^2} \leq M$ for any $k \geq 0$,  then for any $s \in \cS$, 
\begin{align}
 \tsum_{t=1}^k \theta_t \rbr{\hat{\cV}^{\pi_t}(s) - \cV^{\pi_t}(s)}
&  \leq \tsum_{t=1}^k \theta_t 
\hat{\cV}^{\pi_t}(s) - \cV^*(s) \nonumber \\
& \leq \rbr{\tsum_{t=1}^k \beta_t}^{-1}  \tsum_{t=1}^k \frac{\beta_t^2 \gamma^2 \zeta^2  \norm{\hat{\cV}^{\pi_t}}_\infty^2}{2 \mu_w \lambda_{t-1} (1-\gamma)}
+  \rbr{\tsum_{t=1}^k \beta_t}^{-1} \frac{\lambda_k \overline{w}}{1-\gamma} \nonumber \\
& ~~~ + \tsum_{t=1}^k \frac{\theta_t}{1-\gamma}  \EE_{s' \sim d_{s}^{\pi^*}} \sbr{ \delta_t(s', \pi^*(s'))} 
+ \tsum_{t=1}^k \theta_t \rbr{\hat{\cV}^{\pi_t}(s) - \cV^{\pi_t}(s)}. \label{eq_generic_prop_stoch}
\end{align}
\end{lemma}

\begin{proof}
%Let us define $\DD^{\pi(s)}  = \pi(s)$ for any policy $\pi$.
%For notational clarity we will sometimes use these two quantities interchangeably. 
Taking the telescopic sum of \eqref{raw_progress_ineq_pda_stoch}  yields
% \yan{change the annotation to the following, and modify the ensuing results on dependence of $\mu_w$ accordingly}
\begin{align*}
\hat{\Phi}_0 (s, \pi_1(s)) & \leq 
\hat{\Phi}_k(s, \pi_{k+1}(s)) 
- \tsum_{t=1}^k \lambda_{t-1} \cB_s( \pi_t(s), \pi_{t+1}(s)) 
- \tsum_{t=1}^k \beta_t \hat{\phi}_t (s, \pi_{t+1}(s)) \\
& \overset{(a)}{\leq} 
 \hat{\Phi}_k(s, \pi_{k+1}(s)) 
-   \tsum_{t=1}^k \tsum_{a \in \cA} \frac{\lambda_{t-1} \mu_w \vartheta(a|s)}{2} \norm{\DD^{\pi_{t+1}(s)}_a - \DD^{\pi_t(s)}_a }_1^2  \\
& ~~~ - \tsum_{t=1}^k \beta_t  \gamma \zeta \inner{\pi_{t+1}(s) - \pi_t(s)}{\hat{\cV}^{t}_{\vartheta, s}} \\
& \overset{(b)}{=} 
 \hat{\Phi}_k(s, \pi_{k+1}(s)) 
-  \tsum_{t=1}^k \tsum_{a \in \cA}  \frac{\lambda_{t-1} \mu_w \vartheta(a|s)}{2} \norm{\DD^{\pi_{t+1}(s)}_a - \DD^{\pi_t(s)}_a }_1^2  \\
& ~~~  -  \tsum_{t=1}^k \beta_t  \tsum_{a \in \cA}  \gamma \zeta  \vartheta(a|s) \inner{\DD^{\pi_{t+1}(s)}_a - \DD^{\pi_t(s)}_a}{\hat{\cV}^{\pi_t}} \\
& \overset{(c)}{\leq} 
 \hat{\Phi}_k(s, \pi_{k+1}(s)) 
-  \tsum_{t=1}^k \tsum_{a \in \cA}  \frac{\lambda_{t-1} \mu_w \vartheta(a|s)}{2}  \norm{\DD^{\pi_{t+1}(s)}_a - \DD^{\pi_t(s)}_a }_1^2  \\
& ~~~  -   \tsum_{t=1}^k \tsum_{a \in \cA} \beta_t  \gamma \zeta  \vartheta(a|s)   \norm{\hat{\cV}^{\pi_t}}_\infty \norm{\DD^{\pi_{t+1}(s)}_a - \DD^{\pi_t(s)}_a}_1 \\
& \overset{(d)}{\leq} 
 \hat{\Phi}_k(s, \pi_{k+1}(s)) 
+   \tsum_{t=1}^k \tsum_{a \in \cA} \frac{\beta_t^2 \gamma^2 \zeta^2  \norm{\hat{\cV}^{\pi_t}}_\infty^2 \vartheta(a|s)}{2 \mu_w \lambda_{t-1}} \\
& =
 \hat{\Phi}_k(s, \pi_{k+1}(s)) 
+   \tsum_{t=1}^k \frac{\beta_t^2 \gamma^2 \zeta^2  \norm{\hat{\cV}^{\pi_t}}_\infty^2 }{2 \mu_w \lambda_{t-1}} \\
&\overset{(e)}{\leq} \hat{\Phi}_k(s, \DD)  + \tsum_{t=1}^k \frac{\beta_t^2 \gamma^2 \zeta^2  \norm{\hat{\cV}^{\pi_t}}_\infty^2}{2 \mu_w \lambda_{t-1}}, ~ \forall \DD \in \cD_s,
\end{align*}
where $(a)$ follows from the definition of $\hat{\phi}_t$ and \eqref{sc_in_group_norm};
$(b)$ follows from the definition of $\hat{\cV}^t_{\vartheta, s}$;
$(c)$ follows from H\"{o}lder's inequality and the definition of $M$;
$(d)$ follows from Young's inequality; 
and $(e)$ applies \eqref{pda_nhree_point_stoch}.
Since $\beta_0 = 0$ and $w_s(\cdot) \geq 0$, the above inequality implies
\begin{align*}
0 \leq \tsum_{t=1}^k \beta_t {\phi}_t(s, \DD) +  \tsum_{t=1}^k \frac{\beta_t^2 \gamma^2 \zeta^2  \norm{\hat{\cV}^{\pi_t}}_\infty^2}{2 \mu_w \lambda_{t-1}} + \lambda_k w_s(\DD)
+ \tsum_{t=1}^k \beta_t  \delta_t(s,\DD),
\end{align*}
%where $\delta_t \coloneqq \hat{\phi}_n - \phi_t$.
Now consider  aggregating the above inequalities by weights $d_{s}^{\pi^*}$ after taking $\DD=\pi^*(s)$ therein.
From Lemma \ref{lemma_perf_diff} we then obtain  
\begin{align*}
0 \leq \tsum_{t=1}^k \beta_t (1-\gamma) \rbr{\cV^*(s) - \cV^{\pi_t}(s)} 
+  \tsum_{t=1}^k \frac{\beta_t^2 \gamma^2 \zeta^2  \norm{\hat{\cV}^{\pi_t}}_\infty^2}{2 \mu_w \lambda_{t-1}} + \lambda_k \overline{w}
+ \tsum_{t=1}^k \beta_t \EE_{s' \sim d_{s}^{\pi^*}} \sbr{ \delta_t(s', \pi^*(s'))},
\end{align*}
Simple rearrangement of the above inequality yields the desired claim.
%\begin{align*}
%\tsum_{s \in \cS} \rho(s) \tsum_{t=1}^k \theta_t \sbr{
%\cV^{\pi_t}(s) - \cV^*(s)
% }
%& \leq \rbr{\tsum_{t=1}^k \beta_t}^{-1}  \tsum_{t=1}^k \frac{\beta_t^2 \gamma^2 \zeta^2 M^2}{2 \mu_w \lambda_{t-1}}
%+  \rbr{\tsum_{t=1}^k \beta_t}^{-1} \lambda_k \overline{w} \\
%& ~~~ +  \rbr{\tsum_{t=1}^k \beta_t}^{-1}  \rbr{ \tsum_{t=1}^k \beta_t \EE_{s \sim d_{\rho}^{\pi^*}} \sbr{ \delta_t(s, \pi^*(s))} }.
%\end{align*}
\end{proof}

We make the following terminology for accuracy certificate.

\begin{definition}[$\epsilon$-estimator]\label{def_acc_certificate}
For any $\epsilon \geq 0$, we say that a randomized quantity $\hat{\cV}$ is an $\epsilon$-estimator of the robust value function $\cV^*$, defined in \eqref{nature_opt_as_robust_value}, in expectation (resp. in high probability) if 
$-\epsilon \leq \EE \sbr{\hat{\cV}(s)} - \cV^*(s) \leq \epsilon$ (resp. $-\epsilon \leq  {\hat{\cV}}(s) - \cV^*(s) \leq \epsilon$ with high probability) for every $s \in \cS$.
\end{definition}

With Lemma \ref{lemma_generic_prop_stoch} in place, we proceed to establish the convergence of  SFRPE in expectation.

\begin{proposition}\label{thrm_stoch_generic_convergence_expectation}
Fix $\lambda > 0$ and set  
\begin{align}\label{param_choice_stoch_general_expectation}
\beta_k = k^{1/2},   ~ \lambda_k = (k+1)\lambda, ~\forall k \geq 0.
\end{align}
Suppose 
\begin{align}\label{stoch_expecation_conv_bias_condition}
 \norm{\EE_{|k}\hat{\cV}^{\pi_k} - \cV^{\pi_k}}_\infty \leq \varepsilon, ~ \EE_{|k} \sbr{\norm{\hat{\cV}^{\pi_k}}_\infty^2} \leq M, ~ \forall k \geq 1.
\end{align}
Then for any $k \geq 1$, 
\begin{align}\label{ineq_expecatation_sfrpe_opt_gap}
- \varepsilon 
\leq 
\EE \sbr{ \tsum_{t=1}^k \theta_t 
\hat{\cV}^{\pi_t}(s)} - \cV^*(s)  
& \leq  \frac{\gamma^2 \zeta^2 M^2}{\mu_w \lambda \sqrt{k} (1-\gamma)}
+  \frac{4 \lambda  \overline{w}}{\sqrt{k} (1-\gamma)}  + (\frac{2 \gamma \zeta}{1-\gamma} + 1) \varepsilon, ~ \forall s \in \cS.
\end{align}
In particular, taking $\lambda = \frac{\gamma \zeta M}{2 \sqrt{\mu_w \overline{w}}}$ yields 
%\yan{discuss on mis-specify $M$ in choosing $\lambda$ affects the convergence -- looks like linear over-estimate of $M$ leads to linear factor of the upper bound, point out to section on linear function approximation}
\begin{align*}
- \varepsilon 
\leq 
\EE \sbr{ \tsum_{t=1}^k \theta_t 
\hat{\cV}^{\pi_t}(s)} - \cV^*(s)  
& \leq  \frac{4 \gamma \zeta M \sqrt{\overline{w}}}{(1-\gamma) \sqrt{\mu_w k}}  + \frac{3 \varepsilon}{1-\gamma}, ~ \forall s \in \cS.
\end{align*}
%\begin{align*}
%- \varepsilon 
%\leq 
%\EE \sbr{ \tsum_{n=1}^k \theta_n 
%\hat{\cV}^{\pi_n}(s)} - \cV^*(s)  
%& \leq \rbr{\tsum_{n=1}^k \beta_n}^{-1}  \tsum_{n=1}^k \frac{\beta_n^2 \gamma^2 \zeta^2 M^2}{2 \mu_w \lambda_{n-1}}
%+  \rbr{\tsum_{n=1}^k \beta_n}^{-1} \lambda_k \overline{w}  + (2 \gamma \zeta + 1) \varepsilon.
%\end{align*} 
\end{proposition}

\begin{remark}
%In view of the choice of $\lambda$ in Proposition \ref{thrm_stoch_generic_convergence_expectation}, 
When $M$ is unknown, one can instead use an estimate $\hat{M}$ and accordingly choose $\lambda = \frac{\gamma \zeta \hat{M}}{2 \sqrt{\mu_w \overline{w}}}$.
The price for using such an estimate is an $\max \cbr{\hat{M}/M, M/\hat{M}}$ factor increase for the right-hand side of \eqref{ineq_expecatation_sfrpe_opt_gap}.
This would be particularly helpful for our discussion in Section \ref{sec_stoch_linear_approx}, when exact information of $M$ can be difficult to obtain.
\end{remark}

\begin{proof}
Given the definition of $\delta_t(s, \DD)$ in \eqref{def_noise_in_phi}, we obtain 
\begin{align*}
\abs{\EE\sbr{\delta_t(s, \DD)}} 
=\gamma \zeta  \abs{\EE\sbr{ \inner{\DD - \pi_t(s)}{\EE_{|t}\hat{\cV}^{\pi_t}_{\vartheta, s} - {\cV}^{\pi_t}_{\vartheta, s}} }}
\overset{(a)}{\leq} 2 \gamma \zeta \norm{\EE_{|t}\hat{\cV}^{\pi_t} - \cV^{\pi_t}}_\infty
= 2 \gamma \zeta \varepsilon,  ~ \forall s \in \cS, \DD \in \cD_s,
\end{align*}
where $(a)$ follows from a direct application of H\"{o}lder's inequality combined with the definition of $\hat{V}^{\pi}_{\vartheta,s}$ and ${V}^{\pi}_{\vartheta,s}$.
Taking expectation of \eqref{eq_generic_prop_stoch} in Lemma \ref{lemma_generic_prop_stoch} and applying the above inequality yields
\begin{align*}
- \varepsilon 
\leq 
\EE \sbr{ \tsum_{t=1}^k \theta_t 
\hat{\cV}^{\pi_t}(s)} - \cV^*(s)  
& \leq \rbr{\tsum_{t=1}^k \beta_t}^{-1}  \tsum_{t=1}^k \frac{\beta_t^2 \gamma^2 \zeta^2 M^2}{2 \mu_w \lambda_{t-1}(1-\gamma)}
+  \rbr{\tsum_{t=1}^k \beta_t}^{-1} \frac{\lambda_k \overline{w}}{1-\gamma}  + (\frac{2 \gamma \zeta}{1-\gamma} + 1) \varepsilon.
\end{align*} 
The rest of the claims then follow from direct computations after plugging \eqref{param_choice_stoch_general_expectation} into the above inequality. 
%In particular, taking
%\begin{align*}
%\beta_0 = 0, ~ \beta_n = n^{1/2}, ~ \forall n \geq 1; ~ \lambda_n = (n+1)\lambda, ~\forall n \geq 0, 
%\end{align*}
%we obtain 
%\begin{align*}
%- \varepsilon 
%\leq 
%\EE \sbr{ \tsum_{n=1}^k \theta_n 
%\hat{\cV}^{\pi_n}(s)} - \cV^*(s)  
%& \leq  \frac{\gamma^2 \zeta^2 M^2}{\mu_w \lambda \sqrt{k}}
%+  \frac{4 \lambda  \overline{w}}{\sqrt{k}}  + (2 \gamma \zeta + 1) \varepsilon.
%\end{align*}
%Taking $\lambda = \frac{\gamma \zeta M}{2 \sqrt{\mu_w}}$, and using the fact that $\gamma, \zeta \in [0,1]$ yields 
%\begin{align*}
%- \varepsilon 
%\leq 
%\EE \sbr{ \tsum_{n=1}^k \theta_n 
%\hat{\cV}^{\pi_n}(s)} - \cV^*(s)  
%& \leq  \frac{4 \gamma \zeta M}{\sqrt{\mu_w k}}  + 3 \varepsilon.
%\end{align*}
\end{proof}

%Clearly, the above convergence guarantees hinge upon the existence of a policy evaluation operator that certifies condition \eqref{stoch_expecation_conv_bias_condition}.
%We next discuss two policy evaluation operators with this capability,
%and consequently determine the sample complexities of SFRPE when instantiated with these operators. 

Up to now, our discussion is based on the existence of stochastic evaluation operators that can certify noise condition \eqref{stoch_expecation_conv_bias_condition}.
For the remainder of our discussions we proceed to construct such evaluation operators for both tabular setting and with linear function approximation. 
Consequently one can invoke Proposition \ref{thrm_stoch_generic_convergence_expectation} to establish the output of SFRPE being an $\epsilon$-estimator of the robust value function in expectation with $\cO(\abs{\cS} /\epsilon^2)$ sample complexity.
It turns out that one can also obtain a much stronger result with essentially the same number of samples.
In particular, we show  later in this section that the output of SFRPE is an $\epsilon$-estimator of the robust value function with high probability.
This improvement seems to be important for applications of SFRPE in stochastic policy optimization of robust MDPs \cite{li2022first}.

\subsection{Tabular Setting} 

%We then proceed to introduce two policy evaluation operators that can certify condition \eqref{stoch_expecation_conv_bias_condition}.
%Consequently one can invoke Proposition \ref{thrm_stoch_generic_convergence_expectation} to obtain the sample complexity of SFRPE instantiated with these evaluation operators. 
%More importantly, we will also establish direct control over  $ { \tsum_{n=1}^k \theta_n 
%\hat{\cV}^{\pi_n}(s)} - \cV^*(s)  $ with high probability as opposed to the expectation bound in Proposition \ref{thrm_stoch_generic_convergence_expectation}.

%{\bf Simulator-based Evaluation (SE).}
The first evaluation operator presented in Algorithm \ref{alg_se}, named simulator-based evaluation (SE), 
 assumes the access to a so-called simulator of MDP $\cM_{\overline{\PP}}$, and performs an $l$-step process for estimating the value function $\cV^\pi$.
 At each step, SE generates a transition pair $(s,s')$ for each state, where $s'$ denotes the random next state upon committing to an action $a \sim \vartheta(\cdot|s)$ at the state $s$. 
 This transition pairs are then used  to construct auxiliary matrix estimates and update the estimator $\hat{\cV}^\pi$ in an incremental fashion.

\begin{algorithm}[t!]
  \caption{Simulator-based Evaluation (SE)}
  \begin{algorithmic}\label{alg_se}
    \STATE {\bf Initialize:} $R_{-1} = I$ and $\hat{\cV}^{\pi}_{-1} = 0$.
    \STATE Construct $\DD^{\pi, \vartheta} \in \RR^{\abs{\cS} \times \abs{\cS}}$ as
	$
	\DD^{\pi, \vartheta}(s, s') = \tsum_{a \in \cA} \vartheta(a|s) \DD^{\pi(s)}_a (s')    , ~ \forall (s, s') \in \cS \times \cS.
    	$
%	\yan{probably replace this also by sample}
        \FOR{ $i = 0, 1, ... l-1$}
    \STATE Set $\cD = \emptyset$. For each $s \in \cS$, commit action $a \sim \vartheta(\cdot|s)$, and collect $s'$. Save $(s, s')$ into  $\cD$.
    \STATE Construct $\hat{\PP}_i \in \RR^{\abs{\cS} \times \abs{\cS}}$ such that 
    \begin{align*}
     \hat{\PP}_i (s, s') =
    \begin{cases}
     1, ~ (s,s') \in \cD; \\
     0, ~ (s, s') \notin \cD.
     \end{cases}
     \end{align*}
    \STATE Update $R_i = R_{i-1} ((1-\zeta) \hat{\PP}_i + \zeta \DD^{\pi, \vartheta})$, and
    $
    \hat{\cV}^{\pi}_i= \hat{\cV}^{\pi}_{i-1} + \gamma^i R_i \mathfrak{C}.
   $
     \ENDFOR
%    \RETURN $\tsum_{t=1}^k \theta_t \hat{\cV}^{\pi_t}$, where $\theta_t = \beta_t / (\tsum_{t=1}^k \beta_t) $.
\RETURN $\hat{\cV}^{\pi} = \hat{\cV}^{\pi}_l$.
     \end{algorithmic}
\end{algorithm}

%The SE estimator is defined as 
%\begin{align*}
%\hat{\cV}^{\pi}(s) = - \hat{V}^{\vartheta}_{\PP^{\pi}}(s), ~ \text{where} ~ \hat{V}^{\vartheta}_{\PP^{\pi}}(s) = \tsum_{t=0}^{l} \gamma^t c(S_t, A_t).
%\end{align*}
%We make the following immediate observations.

\begin{proposition}\label{prop_se_properties}
For any fixed policy $\pi$, let $\xi$ denote the sample used by the SE operator, then  
\begin{align*}
\norm{\EE_{\xi} \hat{\cV}^{\pi} - \cV^{\pi}}_\infty \leq \frac{\gamma^l}{1-\gamma},~
\norm{\hat{\cV}^{\pi}}_\infty \leq \frac{1}{1-\gamma}.
\end{align*}
\end{proposition}

\begin{proof}
It should be clear that $\EE_{|i} \sbr{ (1-\zeta) \hat{\PP}_i + \zeta \DD^{\pi, \vartheta}} = \mathtt{P}^\pi$,
where $\mathtt{P}^{\pi}$ denotes the transition matrix of the state chain $\cbr{S_t}$ induced by $\vartheta$ within $\cM_{\PP^\pi}$, with $\PP^\pi$ defined as in \eqref{kernel_defined_by_nature_policy}.
Consequently 
from the definition of $R_i$, we obtain
\begin{align*}
\EE_\xi \sbr{R_i}  =  \rbr{ \mathtt{P}^{\pi}}^i  ~ ; ~ \norm{R_i}_\infty  \leq 1, ~ \forall 0 \leq i \leq l,
\end{align*} 
where $\norm{R_i}_\infty$ denotes the matrix $\norm{\cdot}_\infty$ norm.
In addition, one also has 
\begin{align*}
\cV^\pi = (I - \gamma \mathtt{P}^\pi)^{-1}  \mathfrak{C} =  \tsum_{i=0}^\infty \gamma^i \rbr{\mathtt{P}^{\pi}}^i \mathfrak{C}.
\end{align*}
The desired claim then follows immediately from the above two observations and $\norm{\mathfrak{C}}_\infty \leq 1$.
%Consequently, it holds that 
%\begin{align*}
%\norm{\EE_\xi \sbr{\hat{\cV}^{\pi}} - \cV^\pi}_\infty \leq \frac{\gamma^l}{1-\gamma},~
%\norm{\hat{\cV}^{\pi}} \leq \frac{1}{1-\gamma}. 
%\end{align*}
%The rest
%\begin{align*}
%\abs{\EE_\xi \hat{V}^{\vartheta}_{\PP^{\pi}}(s) - V^{\vartheta}_{\PP^{\pi}}(s)} \leq \frac{\gamma^l}{1-\gamma},
%~ \norm{\hat{V}^{\vartheta}_{\PP^{\pi}}}_\infty \leq \frac{1}{1-\gamma}.
%\end{align*}
%The rest of the claims then follows from the above relation.
\end{proof}

\begin{remark}
It should be noted that one can also avoid the construction of $\DD^{\pi, \vartheta}$ via sampling. 
Namely, in addition to the sampled $(s, a, s')$, one also samples $s'' \sim \DD^{\pi(s)}_a(\cdot)$. 
Then $\hat{\DD}^{\pi, \vartheta}_i$ can be constructed in the same way as $\PP_i$ using transition pair $(s, s'')$. 
Accordingly matrix $R_i$ is updated by $R_i = R_{i-1} ((1-\zeta) \hat{\PP}_i + \zeta \hat{\DD}_i^{\pi, \vartheta})$.
\end{remark}

%{\bf Independent Trajectories (IT).}
%IT assumes the access to a so-called simulator. 
%Specifically, for any to be evaluated policy $\pi$, IT generates a trajectory $\xi^{\pi}_{\vartheta}(s)$ of length $l$ as follows: 
%\begin{align*}
%\xi^{\pi}(s) = (S_0=s, A_0, ,  S_1, A_1, , \ldots, S_{l}, A_l ), ~ \text{where} ~ A_t \sim \vartheta(\cdot| S_t), S_{t+1} \sim \PP^{\pi}_{S_t, A_t},
%\end{align*} 
%where $\PP^{\pi}$ is defined as in \eqref{kernel_defined_by_nature_policy}.
%Let us denote $\xi = \cbr{\xi^{\pi}_{\vartheta}(s)}_{s \in \cS}$, then 
%the IT estimator is defined as 
%\begin{align*}
%\hat{\cV}^{\pi}(s) = - \hat{V}^{\vartheta}_{\PP^{\pi}}(s), ~ \text{where} ~ \hat{V}^{\vartheta}_{\PP^{\pi}}(s) = \tsum_{t=0}^{l} \gamma^t c(S_t, A_t).
%\end{align*}
%We make the following immediate observations.
%
%\begin{proposition}\label{prop_se_properties}
%For any fixed policy $\pi$, SE operator yields 
%\begin{align*}
%\norm{\EE_{\xi} \hat{\cV}^{\pi} - \cV^{\pi}}_\infty \leq \frac{\gamma^l}{1-\gamma},~
%\norm{\hat{\cV}^{\pi}}_\infty \leq \frac{1}{1-\gamma}.
%\end{align*}
%\end{proposition}
%
%\begin{proof}
%It is clear that from the definition of $\PP^{\pi}$ and $\xi$, we obtain 
%\begin{align*}
%\abs{\EE_\xi \hat{V}^{\vartheta}_{\PP^{\pi}}(s) - V^{\vartheta}_{\PP^{\pi}}(s)} \leq \frac{\gamma^l}{1-\gamma},
%~ \norm{\hat{V}^{\vartheta}_{\PP^{\pi}}}_\infty \leq \frac{1}{1-\gamma}.
%\end{align*}
%The rest of the claims follow trivially from \eqref{eq_nature_value_as_player_value}.
%\end{proof}

We are now ready to establish the sample complexity of SFRPE with the SE operator. 
As our first result, we establish an $\cO(\abs{\cS}/\epsilon^2)$ sample complexity for SFRPE to output an $\epsilon$-estimator (in expectation) of the robust value function.

\begin{theorem}\label{thrm_sample_se_expectation}
Suppose SFRPE is instantiated with the SE operator, and 
\begin{align*}
\beta_k = k^{1/2},   ~ \lambda_k =  \frac{ (k+1) \gamma \zeta }{2 (1-\gamma) \sqrt{\mu_w \overline{w}}}, ~\forall k \geq 0.
\end{align*}
For any $\epsilon > 0$, to find an approximate robust value such that 
\begin{align*}
-\epsilon \leq \EE \sbr{\tsum_{t=1}^k \theta_t \hat{\cV}^{\pi_t}}(s) - \cV^{*}(s) \leq \epsilon,  ~ s \in \cS,
\end{align*}
SFRPE needs at most $k = 1 + \frac{64 \gamma^2 \zeta^2  \overline{w}}{\mu_w (1-\gamma)^4 \epsilon^2}$ iterations.
The total number of samples can be bounded by 
\begin{align}\label{eq_num_samples_se_expectation}
{\cO} \rbr{ \frac{\gamma^2 \zeta^2 \abs{\cS} \overline{w} \log(1/\epsilon) }{\mu_w \epsilon^2 (1-\gamma)^5 } + \frac{\abs{\cS}}{1-\gamma} \log \rbr{\frac{1}{\epsilon}} }.
\end{align}
In particular, when the distance generating function $w_s(\cdot)$ is set as in \eqref{dgf_negative_entropy}, the number of samples required is bounded by 
\begin{align*}
{\cO} \rbr{ \frac{\gamma^2 \zeta^2 \abs{\cS}  \log \abs{\cS} \log(1/\epsilon) }{ \epsilon^2 (1-\gamma)^5 } + \frac{\abs{\cS}}{1-\gamma} \log \rbr{\frac{1}{\epsilon}} }.
\end{align*}
\end{theorem}

\begin{proof}
Given Lemma \ref{lemma_value_correspondence} and  Proposition \ref{thrm_stoch_generic_convergence_expectation}, it suffices to make sure 
\begin{align*}
\varepsilon \leq \frac{\epsilon(1-\gamma)}{6}, ~ \frac{4 \gamma \zeta M \sqrt{\overline{w}}}{\sqrt{\mu_w k}} \leq \frac{\epsilon(1-\gamma)}{2},
\end{align*}
which in view of Proposition \ref{prop_se_properties}, can be readily satisfied by 
\begin{align*}
l = \frac{1}{1-\gamma} \log \rbr{\frac{6(1-\gamma)^2}{\epsilon}}, ~ k = 1+ \frac{64 \gamma^2 \zeta^2 M^2 \overline{w}}{\mu_w (1-\gamma)^2 \epsilon^2} = 1 + \frac{64 \gamma^2 \zeta^2  \overline{w}}{\mu_w (1-\gamma)^4 \epsilon^2}.
\end{align*}
Consequently, the total number of samples required is bounded by 
\begin{align*}
\abs{\cS} \cdot l \cdot k = {\cO} \rbr{ \frac{\gamma^2 \zeta^2 \overline{w} \abs{\cS} \log(1/\epsilon)}{\mu_w \epsilon^2 (1-\gamma)^5 } + \frac{\abs{\cS}}{1-\gamma} \log \rbr{\frac{1}{\epsilon}} }.
\end{align*}
The rest of the claim follows from \eqref{bregman_divergence_negative_entropy_lb_ub} and \eqref{bregman_divergence_negative_entropy_sc}, from which we conclude the proof.
\end{proof}

%\yan{need some remark on how the sample complexity scales with $\zeta$, and how finite samples are needed despite have continuous action space}
In view of Theorem \ref{thrm_sample_se_expectation}, it is worth noting here that despite SFRPE being applied to solve the MDP $\mathfrak{M}$ of nature with continuous action space,
its sample complexity is independent of the action space.
 This is in sharp contrast when solving general MDPs, where linear dependence on the size of the action space is necessary.
The obtained sample complexity decomposes into two terms that clearly delineates the role of robustness in affecting the sample complexity.
In particular, the first term corresponds to the price we pay for robustness,
and the second term corresponds to the number of samples needed for estimating non-robust value function that is an $\epsilon$-estimator in expectation. 

We next establish the convergence of SFRPE instantiated with the SE operator in high probability. 

\begin{theorem}\label{thrm_stoch_se_high_prob}
Suppose SFRPE is instantiated with the SE operator, and 
\begin{align}\label{se_high_prob_param_choice}
\beta_k = k^{1/2},   ~ \lambda_k =  \frac{ (k+1) \gamma \zeta }{2 (1-\gamma) \sqrt{\mu_w \overline{w}}}, ~\forall k \geq 0.
\end{align}
Then for any $k \geq 0$ and any $\delta \in (0,1)$, with probability at least $1-\delta$ we have 
\begin{align}
& -\sbr{ \frac{\gamma^l}{1-\gamma} + \frac{4}{(1-\gamma)\sqrt{k}} \log \rbr{\frac{4 \abs{\cS} }{\delta}}}  \nonumber \\
 \leq &  \tsum_{t=1}^k \theta_t 
\hat{\cV}^{\pi_t}(s) - \cV^*(s) \nonumber \\
 \leq  & \frac{4 \gamma \zeta  \sqrt{\overline{w}}}{(1-\gamma)^2 \sqrt{\mu_w k}}
+  \frac{\gamma \zeta}{1-\gamma} 
\rbr{
\frac{2 \gamma^l}{1-\gamma} 
+ \frac{8}{(1-\gamma) \sqrt{k}} \sqrt{  \log(\frac{2 \abs{\cS} }{\delta})}
}
+ 
 \frac{\gamma^l}{1-\gamma} + \frac{4}{(1-\gamma)\sqrt{k}} \sqrt{ \log \rbr{\frac{4 \abs{\cS} }{\delta}}}. \label{high_prob_err_bound_se}
\end{align}
Moreover, the total number of samples required by SFRPE to output $-\epsilon \leq \tsum_{t=1}^k \theta_t 
\hat{\cV}^{\pi_t}(s) - \cV^*(s) \leq \epsilon$ with at least probability $1-\delta$  is bounded by 
\begin{align}\label{se_sample_high_prob}
{\cO} \rbr{
\frac{\gamma^2 \zeta^2 \overline{w} \abs{\cS} \log(1/\epsilon) }{(1-\gamma)^5 \mu_w \epsilon^2} \log \rbr{\frac{\abs{\cS}}{\delta}} 
+ \frac{\abs{\cS}  \log(1/\epsilon)}{(1-\gamma)^3 \epsilon^2} \log \rbr{\frac{\abs{\cS}}{\delta}} 
}.
\end{align}
In particular, when the distance generating function $w_s(\cdot)$ is set as in \eqref{dgf_negative_entropy},  the total number of samples required can be bounded by 
\begin{align*}
{\cO} \rbr{
\frac{\gamma^2 \zeta^2  \abs{\cS} \log \abs{\cS}  \log(1/\epsilon)}{(1-\gamma)^5  \epsilon^2} \log \rbr{\frac{\abs{\cS}}{\delta}} 
+ \frac{\abs{\cS}  \log(1/\epsilon)}{(1-\gamma)^3 \epsilon^2} \log \rbr{\frac{\abs{\cS}}{\delta}} 
}.
\end{align*}
\end{theorem}

\begin{proof}
%Let us first consider bounding 
%$\tsum_{t=1}^k \theta_t \rbr{\hat{\cV}^{\pi_t}(s) - \cV^{\pi_t}(s)}$.
%It can be directly verified that 
%\begin{align}\label{eq_opt_error_accumulation}
%\rbr{\tsum_{t=1}^k \beta_n}^{-1}  \tsum_{t=1}^k \frac{\beta_n^2 \gamma^2 \zeta^2  \norm{\hat{\cV}^{\pi_t}}_\infty^2}{2 \mu_w \lambda_{n-1} (1-\gamma)} \leq  \frac{4 \gamma \zeta M \sqrt{\overline{w}}}{(1-\gamma) \sqrt{\mu_w k}}.
%\end{align}
Fixing $s \in \cS$, for any $\delta > 0$, from Proposition \ref{prop_se_properties}, applying Azuma–Hoeffding inequality yields 
%\yan{expand this, skipped definition of $\theta_t$ in the computation, need to follow the second part}
\begin{align*}
\abs{ \tsum_{t=1}^k \theta_t \rbr{\hat{\cV}^{\pi_t}(s) - \cV^{\pi_t}(s)}  } 
 \leq 
\tsum_{t=1}^k  \frac{\theta_t \gamma^l}{1-\gamma} 
+ \frac{1}{1-\gamma} \sqrt{2 \tsum_{t=1}^k \theta_t^2 \log(\frac{4}{\delta})}
 \overset{(a)}{\leq} \frac{\gamma^l}{1-\gamma} + \frac{4}{(1-\gamma)\sqrt{k}} \sqrt{\log \rbr{\frac{4}{\delta}}}, 
\end{align*}
with probability $1 - \delta/ 2$,
where $(a)$ follows from the definition of $\theta_t$ in \eqref{def_theta} together with $\beta_n = n^{1/2}$.
Applying union bound  yields 
\begin{align}\label{se_high_prob_accum_noise_1}
\abs{ \tsum_{t=1}^k \theta_t \rbr{\hat{\cV}^{\pi_t}(s) - \cV^{\pi_t}(s)} }\leq \frac{\gamma^l}{1-\gamma} + \frac{4}{(1-\gamma)\sqrt{k}} \log \rbr{\frac{4 \abs{\cS} }{\delta}}, ~ \forall s \in \cS, 
\end{align}
with probability $1 - \delta/ 2$.
%We proceed to bound $\tsum_{t=1}^k \frac{\theta_t}{1-\gamma}  \EE_{s' \sim d_{s}^{\pi^*}} \sbr{ \delta_t(s', \pi^*(s'))} $.
Fixing $\DD \in \cD_s$, by definition, we have 
\begin{align*}
\delta_t(s,\DD) = \gamma \zeta \tsum_{a \in \cA}  \vartheta(a|s) \inner{\DD_a - \DD^{\pi_t(s)}_a}{\hat{\cV}^{\pi_t} - \cV^{\pi_t}}.
\end{align*}
From Proposition \ref{prop_se_properties}, it is immediate that 
\begin{align*}
\abs{\EE_{|t}  \tsum_{a \in \cA}  \vartheta(a|s)   \inner{\DD_a - \DD^{\pi_t(s)}_a}{\hat{\cV}^{\pi_t} - \cV^{\pi_t}}} \leq \frac{2 \gamma^l}{1-\gamma} ,~
 \abs{  \tsum_{a \in \cA}  \vartheta(a|s) \inner{\DD_a - \DD^{\pi_t(s)}_a}{\hat{\cV}^{\pi_t} - \cV^{\pi_t}}} \leq \frac{2}{1-\gamma}.
\end{align*}
Consequently, applying Azuma–Hoeffding inequality yields
\begin{align*}
\tsum_{t=1}^k   \theta_t   \tsum_{a \in \cA}  \vartheta(a|s) \inner{\DD_a - \DD^{\pi_t(s)}_a}{\hat{\cV}^{\pi_t} - \cV^{\pi_t}}
\leq \frac{2 \gamma^l}{1-\gamma} 
+ \frac{2}{1-\gamma} \sqrt{2 \tsum_{t=1}^k \theta_t^2 \log(\frac{2}{\delta})} 
\leq 
\frac{2 \gamma^l}{1-\gamma} 
+ \frac{8}{1-\gamma} \sqrt{\log(\frac{2}{\delta})} ,
\end{align*}
with probability $1- {\delta}/{2}$, for any $s \in \cS$.
Applying union bound again, we obtain 
\begin{align*}
\tsum_{t=1}^k  \theta_t   \tsum_{a \in \cA}  \vartheta(a|s) \inner{\DD_a - \DD^{\pi_t(s)}_a}{\hat{\cV}^{\pi_t} - \cV^{\pi_t}}
\leq \frac{2 \gamma^l}{1-\gamma} 
+  \frac{8}{1-\gamma} \sqrt{\log(\frac{2 \abs{\cS}}{\delta})} , ~\forall s \in \cS , 
%+ \frac{2}{1-\gamma} \sqrt{2 \tsum_{t=1}^k \theta_t^2 \log(\frac{2 \abs{\cS} \abs{\cA}}{\delta})}, ~ \forall (s,a) \in \cS \times\cA, 
\end{align*}
with probability $1- {\delta}/{2}$.
Setting $\DD = \pi^*(s)$ in the above inequality yields  
\begin{align}
\tsum_{t=1}^k \frac{\theta_t}{1-\gamma}  \EE_{s' \sim d_{s}^{\pi^*}} \sbr{ \delta_t(s', \pi^*(s'))}
%&\leq 
%\frac{\gamma \zeta}{1-\gamma} 
%\rbr{
%\frac{2 \gamma^l}{1-\gamma} 
%+ \frac{2}{1-\gamma} \sqrt{2 \tsum_{t=1}^k \theta_t^2 \log(\frac{2 \abs{\cS} \abs{\cA}}{\delta})}
%}
% \nonumber  \\
 & \leq 
 \frac{\gamma \zeta}{1-\gamma} 
\rbr{
\frac{2 \gamma^l}{1-\gamma} 
+ \frac{8}{(1-\gamma) \sqrt{k}} \sqrt{  \log(\frac{2 \abs{\cS} }{\delta})}
}. \label{se_high_prob_accum_noise_2}
\end{align}
%where the last inequality follows again from $\beta_n = n^{1/2}$ and the definition of $\theta_t$.
Substituting \eqref{se_high_prob_param_choice}, \eqref{se_high_prob_accum_noise_1}, and \eqref{se_high_prob_accum_noise_2}
into \eqref{eq_generic_prop_stoch}, and noting that $\norm{\hat{\cV}^{\pi_t}}_\infty \leq \frac{1}{1-\gamma}$ therein yields \eqref{high_prob_err_bound_se}.
Finally, \eqref{se_sample_high_prob} follows from a similar argument as in Theorem \ref{thrm_sample_se_expectation} and applying \eqref{high_prob_err_bound_se} and Proposition~\ref{prop_se_properties}.
%and using $\beta_n = n^{1/2}$, $\lambda_n = (n+1) \lambda$ with $\lambda = \frac{\gamma \zeta M}{2 \sqrt{\mu_w \overline{w}}}$ and $M = \frac{1}{1-\gamma}$,
%we obtain
%\begin{align*}
%& -\sbr{ \frac{\gamma^l}{1-\gamma} + \frac{4}{(1-\gamma)\sqrt{k}} \log \rbr{\frac{4 \abs{\cS} }{\delta}}} \\
% \leq &  \tsum_{t=1}^k \theta_t 
%\hat{\cV}^{\pi_t}(s) - \cV^*(s) \nonumber \\
% \leq  & \frac{4 \gamma \zeta  \sqrt{\overline{w}}}{(1-\gamma)^2 \sqrt{\mu_w k}}
%+  \frac{\gamma \zeta}{1-\gamma} 
%\rbr{
%\frac{2 \gamma^l}{1-\gamma} 
%+ \frac{8}{(1-\gamma) \sqrt{k}} \sqrt{  \log(\frac{2 \abs{\cS} \abs{\cA}}{\delta})}
%}
%+ 
% \frac{\gamma^l}{1-\gamma} + \frac{4}{(1-\gamma)\sqrt{k}} \log \rbr{\frac{4 \abs{\cS} }{\delta}}.
%\end{align*}
%The proof is then completed.
\end{proof}

A few remarks are in order for interpreting Theorem \ref{thrm_stoch_se_high_prob}.
First, it is clear that the sample complexity in \eqref{se_sample_high_prob} is comparable to that of \eqref{eq_num_samples_se_expectation}, while the accuracy certificate is now stated with high probability instead of expectation. 
Second, similar to Theorem \ref{thrm_sample_se_expectation}, the sample complexity in \eqref{se_sample_high_prob} consists of two terms of different roots.
The first term corresponds to the price of robustness,
and the second term corresponds to the number of samples required for estimating the standard value function up to $\epsilon$-accuracy with high probability. 
Clearly, when $\zeta =0$ the established sample complexities are tight for evaluating the standard value function in both expectation and in high probability.
It is also important to note that 
in view of \eqref{se_sample_high_prob},  when $\zeta = \cO (1-\gamma)$,  the robust value function can be evaluated with the same number of samples as for evaluating standard value function. 
Consequently there is no additional price of robustness for robust policy evaluation with small-sized ambiguity sets.

The development of SFRPE in this section appears to be new in several aspects. 
The convergence of SFRPE  is established in both expectation and high probability, while existing development of stochastic robust policy evaluation only establish convergence in expectation for mean-squared error \cite{li2022first}. 
Additionally, the impact of the ambiguity set size on the sample complexity has not been previously reported.
%In addition, the size of the ambiguity set in affecting the sample complexity has not been reported before.  
Though these results already hint upon potential benefits of the SFRPE  framework, in the next section we proceed to demonstrate its true advantage of scaling robust policy evaluation to large-scale problems, 
a scenario that appears yet to have an algorithmic solution.
%where no previous algorithmic solution exists.

%
%\yan{need some remark on how the sample complexity scales with $\zeta$,
%\begin{itemize}
%\item the first term of the sample complexity corresponds to the price we pay for robustness -- notably when $\zeta = \cO(1-\gamma)$ there is no price to pay for robustness!
%\item the second term corresponds to the one for standard evaluation 
%\end{itemize}
%}
%
%
%\yan{comparison to prior work: this is the first work where linear dependence on state space size can be obtained. previous approach, aside from making strong assumptions, is inherently subject to the exploration issue over state space}
%
%
%

%!TEX root = ./robust_pe.tex

\subsection{SFRPE with Linear Function Approximation}\label{sec_stoch_linear_approx}

Unless stated otherwise, going forward we let $\norm{\cdot} = \norm{\cdot}_2$.
For robust MDP with large state space, exact policy evaluation becomes prohibitive. 
In this case one can instead seek to learn a linearly parameterized $\cV^{\pi}_{\theta}(\cdot) \coloneqq   \psi(\cdot)^\top \theta$ that approximate $\cV^{\pi}(\cdot)$ well, 
where $\psi: \cS \to \RR^d$ is the so-called feature map, and we assume without loss of generality that $\norm{\psi(\cdot)}\leq 1$.
In view of Lemma \ref{lemma_value_correspondence}, this is equivalent to approximation of $V^{\vartheta}_{\PP^\pi}$ by $- \psi(\cdot)^\top \theta$. 
As $V^{\vartheta}_{\PP^\pi}$ itself is the value function of $\vartheta$ within MDP $\cM_{\PP^\pi}$, we consider solving\footnote{
Interested readers might suggest directly learning the robust value function of $V^{\vartheta}_r$ -- the ultimate goal of this manuscript -- by simply formulating a least-squares objective that fits the mean-squared robust  Bellman error. 
Though intuitively appealing, this perspective can be computationally intractable. Namely, one can easily construct a $\mathrm{s}$-rectangular robust MDP instance where the resulting  least-squares objective of robust mean-squared Bellman error  is non-convex even in the tabular setting.
Accordingly one can only seek to find the stationary point of the least-square objective  \cite{roy2017reinforcement}.
}: 
\begin{align}\label{ls_objective}
\textstyle
\min_{\theta \in \RR^d} \cbr{g(\theta) \coloneqq  \norm{\Psi \theta - \gamma \mathtt{P}^\pi \Psi \theta - \mathfrak{C}}_{\nu}^2 },
\end{align}
which corresponds to the mean-squared Bellman error of $\vartheta$ within MDP $\cM_{\PP^\pi}$.
Here $\Psi \in \RR^{\abs{\cS} \times d}$ denotes the feature matrix induced by $\psi$ applied to every state, and $\mathtt{P}^{\pi}$ denotes the  transition matrix of the state chain $\cbr{S_t}$ induced by $\vartheta$ within $\cM_{\PP^\pi}$.
It is important to note that here we do not know $\mathtt{P}^\pi$.
Instead, we only assume the access to sample from $\cM_{\overline{\PP}}$.
We introduce the following standard assumption on the $\nu$ and $\psi$.

\begin{algorithm}[t]
  \caption{Stochastic Least-squares Policy Evaluation (SLPE)}
  \begin{algorithmic}\label{alg_lspe}
    \STATE {\bf Input:} Stepsizes $\cbr{\eta_t}$.
    \STATE {\bf Initialize:} $\theta_0 \in \RR^d$.
        \FOR{ $t = 0, 1, ... T-1$}
    \STATE Sample $s_t \sim \nu$, commit action $a_t \sim \vartheta(\cdot |s_t)$. Sample independent $x_t, x_t' \sim \overline{\PP}_{s_t, a_t}$
    and $y_t, y_t' \sim \DD^{\pi(s_t)}_{a_t}$. 
    \STATE Update:
    \begin{align*}
    \theta_{t+1} = \theta_t - \eta_t \sbr{
    \psi(s_t)^\top \theta_t - \gamma \rbr{(1-\zeta) \psi(x_t') + \zeta \psi(y_t') }^\top \theta_t - \mathfrak{C}(s_t)
    }
    \sbr{
    \psi(s_t) - \gamma \rbr{ (1-\zeta) \psi(x_t) +  \zeta \psi(y_t) }
    }
    \end{align*}
     \ENDFOR
%    \RETURN $\tsum_{t=1}^k \theta_t \cV^{\pi_t}$, where $\theta_t = \beta_t / (\tsum_{t=1}^k \beta_t) $.
\RETURN $\theta_T$, and $\hat{\cV}^{\pi}(\cdot) = \cV^{\pi}_{\theta_T}(\cdot) \coloneqq  \psi(\cdot)^\top \theta_T$.
     \end{algorithmic}
\end{algorithm}

%\yan{should mention that directly using robust bellman operator to form least square objective leads to non-convex objective}

\begin{assumption}\label{assump_sampling_and_feature}
The distribution $\nu$ has full support and the feature matrix $\Psi \in \RR^{\abs{\cS} \times d}$ is non-singular.
That is, $\diag(\nu) \succ 0$ and $\sigma_{\min}(\Psi) > 0$.
\end{assumption}

Following Assumption \ref{assump_sampling_and_feature},
it holds that $\mu = \lambda_{\min} ( \Psi^\top (I - \gamma \mathtt{P}^{\pi})^\top \diag(\nu) (I - \gamma \mathtt{P}^{\pi}) \Psi) ) > 0$.
We also denote $L = \lambda_{\max} \rbr{ \Psi^\top (I - \gamma \mathtt{P}^{\pi})^\top \diag(\nu) (I - \gamma \mathtt{P}^{\pi}) \Psi}$.
Let us denote $\theta^\pi$ as the unique solution of \eqref{ls_objective}, and use 
\begin{align*}
\textstyle
\varepsilon_{\mathrm{approx}} \coloneqq \sup_{\pi \in \Pi} \norm{ \cV^{\pi}_{\theta^\pi} - \cV^{\pi} }_\infty
%= \sup_{\pi \in \Pi} \norm{ V^{\vartheta}_{\PP^\pi} - \Psi \theta^{\pi} }_\infty
\end{align*}
 to characterize the function approximation error of $\cV^{\pi}$. 
Clearly, when $\Psi = I_{\abs{\cS}}$, we have $\varepsilon_{\mathrm{approx}} = 0$.
Our ensuing discussions will often make use of the following quantities to simplify presentation:
\begin{align*}
\textstyle
%c_1 = 4 \norm{\theta^{\pi}} + 2, 
r_{\Theta} = \max_{\pi \in \Pi} \norm{\theta^\pi} , 
~
c_1 = 4 r_{\Theta} + 2 .
\end{align*}
With Assumption \ref{assump_sampling_and_feature} it holds that $r_{\Theta} < \infty$ as $\theta^{\pi}$ is a continuous mapping from $\Pi$ to $\RR^d$.

The stochastic least-squares policy evaluation (SLPE) method (Algorithm \ref{alg_lspe}) can be viewed as solving \eqref{ls_objective} by stochastic gradient descent. 
It can be seen that SLPE requires a simulator of $\cM_{\PP^\pi}$ to draw sample $s_t$, a condition that mainly serves to keep the technical discussion concise. With a slightly more complex analysis  one can also sample $s_t$ by following the trajectory of $\vartheta$ within $\cM_{\PP^\pi}$  \cite{kotsalis2020simple}.

\begin{lemma}
Define operator $F: \RR^d \to \RR^d$ as 
\begin{align*}
F(\theta) = \rbr{(I - \gamma \mathtt{P}^{\pi}) \Psi}^\top \diag(\nu) (\Psi \theta - \gamma \mathtt{P}^{\pi} \Psi \theta - \mathfrak{C}).
\end{align*}
%Suppose $\diag(\nu) \succ 0$ and $\sigma_{\min}(\Psi) > 0$, then  
Then   $F(\theta^{\pi}) = 0$, and
\begin{align}
\inner{F(\theta)}{\theta - \theta^{\pi}}   \geq  \mu \norm{\theta - \theta^{\pi}}^2. \label{ctd_monotone}
\end{align}
%where  $F(\theta^{\pi}) = 0$.
%Given Assumption \ref{assump_ergodic} we also have $\mu > 0$, $M \succ \mathbf{0}$ and $\theta^{\pi}$ being unique.
\end{lemma}
\begin{proof}
The first part of the claim directly follows from the optimality condition of \eqref{ls_objective}.
In addition, 
\begin{align}
\inner{F(\theta)}{\theta - \theta^{\pi}}
& = \inner{F(\theta) - F(\theta^{\pi})}{\theta - \theta^{\pi}} \nonumber \\
& = \inner{\Psi^\top (I - \gamma \mathtt{P}^{\pi})^\top \diag(\nu)(I-\gamma \mathtt{P}^\pi) \Psi (\theta - \theta^{\pi})}{\theta - \theta^{\pi}} \nonumber \\
& \geq  \mu \norm{\theta - \theta^{\pi}}^2, \nonumber
\end{align}
where the last inequality follows from the definition of $\mu$.
\end{proof}

We proceed by characterizing each step of SLPE and establish the boundedness of iterates in expectation. 

\begin{lemma}\label{lemma_ctd_recursion}
%Define $c_1 = 4 \norm{\theta^{\pi}} + 2$, then
We have
%\begin{align*}
%L= \sigma_{\max}\rbr{\Psi^\top (I - \gamma \mathtt{P}^\pi) \Psi} , ~ & \kappa_2  = \sigma_{\max} \rbr{\Psi^\top (I-\gamma \mathtt{P}^\pi)},  
%\kappa_3 =  \lambda_{\max} \rbr{\Psi^\top \Psi} + 1, \\
%c_1 = 2 \norm{\theta^{\pi}} + 1, 
%~ & \varepsilon_{\mathrm{approx}} = \sup_{\pi \in \Pi} \norm{\Psi \theta^{\pi} - V^{\vartheta}_{\PP^{\pi}}},
%\end{align*}
%where $V^{\vartheta}_{{\PP}^\pi}$ denotes the value function of $\vartheta$ within $\cM_{\mathtt{P}^\pi}$.
\begin{align}
 \EE \sbr{\norm{\theta_{t+1} - \theta^{\pi}}^2 }\leq 
 \rbr{
1 - 2 \eta_t \mu + 32 \eta_t^2  
}
\EE \sbr{\norm{ \theta_t - \theta^{\pi}}^2} + 2 \eta_t^2 c_1^2 .
\label{convergence_mse}
\end{align}
In particular, setting $\eta_t = \eta \leq \frac{\mu}{32}$ yields
%\begin{align}\label{param_for_norm_bound}
%\eta_t = \eta \leq \frac{\mu}{32}.
%%~ (L + \kappa_2) m \rho^\tau \leq \frac{\mu}{4}.
%%T \geq  \frac{1}{\eta \mu} \log 2.
%\end{align}
%Then we have 
\begin{align}\label{mse_bound_each_epoch_iter}
\EE \sbr{\norm{{\theta}_t - \theta^{\pi}}^2} \leq \norm{\theta_0 - \theta^{\pi}}^2 + 
c_1^2.
%+ \frac{2m \kappa_2 \rho^\tau}{\mu} \varepsilon^2_{\mathrm{approx}}.
\end{align}
 
\end{lemma}

\begin{proof}
Clearly,  each update in SLPE (Algorithm \ref{alg_lspe}) is equivalent to the following
\begin{align}\label{ctd_update_equiv_form}
\textstyle
\theta_{t+1} = \argmin_{\theta \in \RR^d} \eta_t \inner{\hat{F}_t(\theta_t)}{\theta} - \frac{1}{2} \norm{\theta - \theta_t}^2,
\end{align}
where 
$
\hat{F}_t (\theta) =  \sbr{
    \psi(s_t)^\top \theta - \gamma \rbr{(1-\zeta) \psi(x_t') + \zeta \psi(y_t') }^\top \theta - \mathfrak{C}(s_t)
    }
    \sbr{
    \psi(s_t) - \gamma \rbr{ (1-\zeta) \psi(x_t) +  \zeta \psi(y_t) }
    }
$.
From the optimality condition of \eqref{ctd_update_equiv_form}, we obtain 
\begin{align}\label{ctd_three_point}
\eta_t \inner{\hat{F}_t(\theta_t)}{ \theta_t - \theta} + \eta_t \inner{\hat{F}_t (\theta_t)}{\theta_{t+1} - \theta_t}
+ \frac{1}{2} \norm{\theta_{t+1} - \theta_t}^2 \leq \frac{1}{2} \norm{\theta - \theta_t}^2 - \frac{1}{2} \norm{\theta - \theta_{t+1}}^2.
\end{align}
We now make the following two observations. First,  from the definition of $(s_t, x_t, x_t', y_t, y_t')$ it is clear that 
\begin{align*}
 \EE_{|t} \sbr{ \hat{F}_t(\theta) }
 = \EE_{|s_t, t} \big[ \sbr{\psi(s_t)^\top \theta - \gamma  \mathtt{P}^\pi_{s_t, \cdot} \Psi \theta - \mathfrak{C}(s_t)}
 \sbr{\psi(s_t) - \gamma \mathtt{P}^\pi_{s_t, \cdot} \Psi } \big]
 = F(\theta), ~ \forall \theta \in \RR^d,
\end{align*}
where $\EE_{|t} \sbr{\cdot}$ denotes the conditional expectation with respect to the $\sigma$-algebra up to (excluding) iteration $t$.
Consequently, by denoting $\delta_t = \hat{F}_t(\theta_t) - F(\theta_t)$, then
\begin{align}
\EE_{|t} \sbr{\delta_t} = 0.
 \label{bias_conditional_expectation}
\end{align}
%where $\EE_{|t} \sbr{\cdot}$ denotes the conditional expectation with respect to the $\sigma$-algebra up to (excluding) iteration $t$ of epoch $e$.
%where 
%\begin{align*}
%L= \sigma_{\max}\rbr{\Psi^\top (I - \gamma \mathtt{P}^\pi) \Psi} , ~ \kappa_2 = \sigma_{\max} \rbr{\Psi^\top (I-\gamma \mathtt{P}^\pi)},  ~ \varepsilon_{\mathrm{approx}} = \norm{\Psi \theta^{\pi} - V^{\vartheta}_{{\PP}^\pi}}.
%\end{align*}
Second, we have 
$ \hat{F}_t(\theta_t)  = \hat{F}_t(\theta_t) - \hat{F}_t(\theta^{\pi}) + \hat{F}_t(\theta^{\pi})  $,
and 
\begin{align}
\norm{ \hat{F}_t(\theta_t) - \hat{F}_t(\theta^{\pi})}
& = \norm{
\sbr{\psi(s_t) -  \gamma \rbr{(1-\zeta) \psi(x_t') + \zeta\psi(y_t')}}^\top  (\theta_t - \theta^\pi)  \sbr{
    \psi(s_t) - \gamma \rbr{ (1-\zeta) \psi(x_t) +  \zeta \psi(y_t) }
    }
} \nonumber  \\
& \leq 4 \norm{\theta_t - \theta^\pi},  \label{norm_bound_diff_stoch_op} \\
\norm{\hat{F}_t(\theta^{\pi})}
& = \norm{
\sbr{
    \psi(s_t)^\top \theta^\pi - \gamma \rbr{(1-\zeta) \psi(x_t') + \zeta \psi(y_t') }^\top \theta^\pi - \mathfrak{C}(s_t)
    }
    \sbr{
    \psi(s_t) - \gamma \rbr{ (1-\zeta) \psi(x_t) +  \zeta \psi(y_t) }
    }
} \nonumber \\
& \leq 4 \norm{\theta^\pi} + 2. \label{norm_bound_stoch_op}
\end{align}
Hence it follows that 
\begin{align*}
\inner{\hat{F}_t(\theta_t)}{\theta_{t+1} - \theta_t} 
& = \inner{  \hat{F}_t(\theta_t) - \hat{F}_t(\theta^{\pi}) + \hat{F}_t(\theta^{\pi})}{\theta_{t+1} - \theta_t} \nonumber \\
& \geq 
- 4 \norm{\theta_t - \theta^{\pi}} \norm{\theta_{t+1} - \theta_t} 
- \norm{\hat{F}_t(\theta^{\pi})} \norm{\theta_{t+1} - \theta_t} \nonumber  \\
& \geq 
-( 4 \norm{\theta_t - \theta^{\pi}} + c_1) \norm{\theta_{t+1} - \theta_t} . 
%\label{ctd_subgrad_inner}
\end{align*}
%where 
%\begin{align*}
%\kappa_3 =  \lambda_{\max} \rbr{\Psi^\top \Psi} + 1,~
%c_1 = 2 \norm{\theta^{\pi}} + 1.
%\end{align*}
Substituting  the above relation into \eqref{ctd_three_point} yields 
\begin{align*}
\eta_t \inner{\hat{F}_t(\theta_t)}{ \theta_t - \theta} - \eta_t ( 4 \norm{\theta_t - \theta^{\pi}} + c_1) \norm{\theta_{t+1} - \theta_t}
+  \frac{1}{2} \norm{\theta_{t+1} - \theta_t}^2 \leq \frac{1}{2} \norm{\theta - \theta_t}^2 - \frac{1}{2} \norm{\theta - \theta_{t+1}}^2,
\end{align*}
which after applying Cauchy-Schwarz inequality and the definition of $\delta_t$, gives 
\begin{align*}
\eta_t \inner{{F}(\theta_t)}{ \theta_t - \theta} - 16 \eta_t^2  \norm{\theta_t - \theta^{\pi}}^2
- \eta_t^2 c_1^2
 \leq \frac{1}{2} \norm{\theta - \theta_t}^2 - \frac{1}{2} \norm{\theta - \theta_{t+1}}^2
 - \eta_t \inner{\delta_t}{\theta_t - \theta}.
\end{align*}
 Setting $\theta = \theta^{\pi}$, and further plugging \eqref{ctd_monotone} into the above relation yields 
\begin{align}\label{norm_recursion_with_noise}
\frac{1}{2} \norm{\theta_{t+1} - \theta^{\pi}}^2 
\leq (\frac{1}{2} - \eta_t \mu + 16 \eta_t^2 ) \norm{\theta_t - \theta^{\pi}}^2
+ \eta_t^2 c_1^2 +  \eta_t \inner{\delta_t}{\theta_t - \theta^{\pi}}.
\end{align}
Taking expectation on both sides and applying  \eqref{bias_conditional_expectation}, we obtain \eqref{ctd_update_equiv_form}.
 Finally, recursive application of \eqref{convergence_mse}  with $\eta_t = \eta \leq \frac{\mu}{32}$ yields 
\begin{align*}
%\label{convergence_mse_clean}
\EE \sbr{\norm{\theta_t - \theta^{\pi}}^2 }
\leq (1-\eta \mu)^t  \EE \sbr{ \norm{\theta_{0} - \theta^{\pi}}^2 } + \frac{2 \eta c_1^2}{\mu} ,
%+ \frac{2m \kappa_2 \rho^\tau}{\mu} \varepsilon^2_{\mathrm{approx}},
\end{align*}
from which the desired claim follows.
%\begin{align*}
%\frac{1}{2} \EE \sbr{\norm{\theta_{t+1} - \theta^{\pi}}^2 }\leq 
% \rbr{
%\frac{1}{2} - \eta_t \mu + 16 \eta_t^2   
%}
%\EE \sbr{\norm{ \theta_t - \theta^{\pi}}^2} + \eta_t^2 c_1^2.
%%\label{convergence_mse}
%\end{align*}
%Further dividing both sides by $\Gamma_{t+1}$, where 
%\begin{align*}
%\Gamma_t = \begin{cases}
%1, ~ & t = 0; \\
%\Gamma_{t-1}  \rbr{1 - 2 \eta_t \mu + 2 \eta_t^2 \kappa_3^2 + 2 \eta_t m (L + \kappa_2) \rho^\tau }, & t \geq 1,
%\end{cases}
%\end{align*}
%and taking the telescopic sum of the resulting inequality, we obtain 
%\begin{align*}
%\frac{1}{2 \Gamma_k} \EE \sbr{ \norm{\theta_k - \theta^{\pi}}^2 } \leq \frac{1}{2} \norm{\theta_0 - \theta^{\pi}}^2 + \tsum_{t=0}^{k-1} \frac{\eta_t^2}{\Gamma_{t+1}} c_1^2
%+ \tsum_{t=0}^{k-1} \frac{\eta_t}{\Gamma_{t+1}} \kappa_2 \rho^\tau \varepsilon^2_{\mathrm{approx}}.
%\end{align*}
%Hence we conclude that 
%\begin{align*}
%\EE \sbr{ \norm{\theta_k - \theta^{\pi}}^2 }
% \leq \Gamma_k  \norm{\theta_0 - \theta^{\pi}}^2
% + 2 \Gamma_k c_1^2 \tsum_{t=0}^{k-1} \frac{\eta_t^2}{\Gamma_{t+1}} 
% + 2 \Gamma_k  \kappa_2 m \rho^\tau \varepsilon^2_{\mathrm{approx}} \tsum_{t=0}^{k-1} \frac{\eta_t}{\Gamma_{t+1}} .
%\end{align*}
\end{proof}

The next lemma establishes the fast bias reduction  of the estimated value function.

\begin{lemma}\label{lemma_bias_linear}
%Let $L = \lambda_{\max} \rbr{ \Psi^\top (I - \gamma \mathtt{P}^{\pi})^\top \diag(\nu) (I - \gamma \mathtt{P}^{\pi}) \Psi}$, and
Denote $\overline{\theta}_t = \EE \sbr{\theta_t}$, and let
\begin{align}\label{param_choice_bias_reduction}
\eta_t = \eta \leq \frac{\mu}{ L^2}.
%T \geq  \frac{1}{\eta \mu} \log 2.
\end{align}
Then we have 
\begin{align}\label{bound_on_parameter_bias}
\norm{\overline{\theta}_t - \theta^{\pi}}^2
\leq 
(1 - \eta \mu)^t \norm{\theta_0 - \theta^{\pi}}^2.
\end{align}
Consequently, 
%by denoting $\cV^{\pi}_{\theta} = - \Psi \theta$, 
we obtain 
\begin{align}\label{bound_on_value_bias}
\norm{\EE \sbr{ \cV^{\pi}_{\theta_t}} - \cV^{\pi}}_\infty \leq  (1 - \eta \mu)^{t/2} \norm{\theta_0 - \theta^{\pi}}
+  \varepsilon_{\mathrm{approx}} .
\end{align}
%\yan{need a statement on the err of the estimated value}
\end{lemma}

\begin{proof}
%Define $\overline{\theta}_t = \EE\sbr{\theta_t}$ and $\overline{\theta} = \EE\sbr{\theta}$. 
Since we have $\theta_{t+1} = \theta_t - \eta_t \hat{F}_t(\theta_t)$, taking expectation on both sides yields 
\begin{align*}
\overline{\theta}_{t+1} = \EE\sbr{\theta_{t+1}}
= \EE\sbr{\theta_t} - \eta_t \EE [\hat{F}_t(\theta_t)] 
& = \EE\sbr{\theta_t}
-  \eta_t \EE [ \EE_{|t} [\hat{F}_t(\theta_t)] ] \\
& = 
 \EE\sbr{\theta_t}
- \eta_t \EE [ F(\theta_t) + \EE_{|t} [\delta_t] ]  \\
& = 
\overline{\theta}_t 
- \eta_t  F(\overline{\theta}_t)
\end{align*}
where the last equality follows from the linearity of $F(\cdot)$ and \eqref{bias_conditional_expectation}.
Consequently, we have 
\begin{align}
\norm{\overline{\theta}_{t+1} - \theta^{\pi}} 
& = \norm{ \overline{\theta}_t 
-  \eta_t  F(\overline{\theta}_t)
- \theta^{\pi}
}  \nonumber \\
& \leq
\norm{\overline{\theta}_t - \theta^{\pi}}^2 - 2 \eta_t \inner{F(\overline{\theta}_t) }{\overline{\theta}_t - \theta^{\pi}}
+   \eta_t^2 \norm{F(\overline{\theta}_t)}^2 . 
\label{bias_progress_raw}
\end{align}
%We proceed to bound $\norm{\EE\sbr{\delta_t}}$ and $\norm{F(\overline{\theta}_t)}$ separately. 
%For $\norm{\EE\sbr{\delta_t}}$, we have 
%\begin{align}
%\norm{\EE[\delta_t]}
%= \norm{\EE[ \EE_{|t}[\delta_t] ]}
%&\overset{(a)}{ =} \lVert
%\EE \big[
% \Psi^\top (\EE_{|t}[\hat{M}] - M)(I -\gamma \mathtt{P}^\pi) \Psi (\theta_t  -  \theta^{\pi} )
% \nonumber \\
%& ~~~~~~ + \Psi^\top (\EE_{|t}[\hat{M}] - M) (I -\gamma \mathtt{P}^\pi)  [\Psi \theta^{\pi} - V^{\vartheta}_{{\PP}^\pi}] 
%\big]
%\rVert \nonumber \\
%& \leq 
%m \rho^\tau \sbr{ L \EE \sbr{\norm{\theta_t - \theta^{\pi}}}
%+  \kappa_2 \varepsilon_{\mathrm{approx}} }, \label{norm_bound_on_bias},
%\end{align}
%where $(a)$ follows from \eqref{bias_conditional_expectation}.
In addition,  it holds that 
\begin{align}
\norm{F(\overline{\theta}_t)}
& = \norm{F(\overline{\theta}_t) - F(\theta^\pi)}  = \norm{\Psi^\top (I - \gamma \mathtt{P}^{\pi})^\top \diag(\nu) (I - \gamma \mathtt{P}^{\pi}) \Psi (\overline{\theta}_t - \theta^{\pi}}  \leq L  \norm{\overline{\theta}_t - \theta^{\pi}}. 
\label{bound_on_op_norm}
\end{align}
Hence by plugging \eqref{bound_on_op_norm} into \eqref{bias_progress_raw}, and applying \eqref{ctd_monotone}, we obtain 
\begin{align}
\norm{\overline{\theta}_{t+1} - \theta^{\pi}}^2 
& \leq (1- 2\eta_t \mu + L^2 \eta_t^2) \norm{\overline{\theta}_t - \theta^{\pi}}^2 
\label{convergence_bias}
\end{align}
The above relation simplifies to \eqref{bound_on_parameter_bias} with choice of $\cbr{\eta_t}$ in \eqref{param_choice_bias_reduction}.
%\begin{align*}
%\norm{\overline{\theta}_t - \theta^{\pi}}^2
%\leq 
%(1 - \eta \mu)^t \norm{\theta_0 - \theta^{\pi}}^2
%+  \varepsilon^2_{\mathrm{approx}} .
%\end{align*}
%Consequently, with the definition of $\overline{\theta}$ and the choice of $T$, we obtain   
%\begin{align*}
%\norm{\overline{\theta}^{(e+1)} - \theta^{\pi}}^2 \leq \frac{1}{2} \norm{\overline{\theta} - \theta^{\pi}}^2 
%+ \varepsilon^2_{\mathrm{approx}} ,
%\end{align*}
%from which \eqref{bound_on_parameter_bias} follows.
%It remains to note that \eqref{param_for_bias_bound_1} and \eqref{param_for_bias_bound_2} can be satisfied by the parameter choice specified in \eqref{param_choice_bias_reduction}.
%, from which we obtain \eqref{bound_on_parameter_bias}.
Finally, \eqref{bound_on_value_bias} follows from the direct application of \eqref{bound_on_parameter_bias},
by noting that 
\begin{align*}
\abs{\EE \sbr{ \cV^{\pi}_{\theta_t}}(s) - \cV^{\pi}(s)}
& \leq 
 \abs{\EE \sbr{ \cV^{\pi}_{\theta_t}}(s) - \cV^\pi_{\theta^\pi}(s)} +  \abs{ \cV^\pi_{\theta^\pi}(s)- \cV^{\pi}(s)}
\\ & \leq 
 \abs{ \rbr{\overline{\theta}_t - \theta^\pi}^\top \psi(s)} +  \abs{ \cV^\pi_{\theta^\pi}(s)- \cV^{\pi}(s)} \\
 & \leq 
  (1 - \eta \mu)^{t/2} \norm{\theta_0 - \theta^{\pi}} +  \varepsilon_{\mathrm{approx}} .
\end{align*}
The proof is then completed.
\end{proof}

With Lemma \ref{lemma_ctd_recursion} and \ref{lemma_bias_linear} in place, we are ready to establish the sample complexity of SFRPE with the SLPE operator, 
in order to output an $\epsilon$-estimator of the robust value function in expectation.
%For notational simplicity, going forward we denote 
%$r_{\Theta} = \max_{\pi \in \Pi} \norm{\theta^\pi} < \infty$.

\begin{theorem}\label{thrm_lspe_expectation}
For any $\epsilon \geq \frac{8 \varepsilon_{\mathrm{approx}} }{(1-\gamma)}$,
let SFRPE be instantiated with SLPE operator with  evaluation parameters
\begin{align*}
\eta \leq \frac{\mu}{L^2 +  32} , ~ \theta_0 = 0, ~ T = \cO \rbr{ \frac{1}{\eta \mu} \log \rbr{\frac{r_\Theta }{ \epsilon}} }, 
\end{align*}
 and  optimization parameters
\begin{align*}
\beta_k = k^{1/2},   ~ \lambda_k = \frac{ (k+1) M \gamma \zeta }{2  \sqrt{\mu \overline{w}}}, ~\forall k \geq 0, 
\end{align*}
where $M^2 \geq 4 \rbr{r_\Theta^2 + c_1^2  + \varepsilon^2_{\mathrm{approx}}} + {2}/{(1-\gamma)^2}$.
To find an approximate robust value such that 
\begin{align*}
-\epsilon \leq \EE \sbr{\tsum_{t=1}^k \theta_t \cV^{\pi_t}}(s) - \cV^*(s) \leq \epsilon,  ~ s \in \cS,
\end{align*}
SFRPE needs at most $k = 1 + \frac{256 \gamma^2 \zeta^2 M^2 \overline{w}}{(1-\gamma)^2 \mu \epsilon^2}$ iterations.
%where $M^2 = 4 \rbr{r_\Theta^2 + 1  + \varepsilon^2_{\mathrm{approx}}} + {2}/{(1-\gamma)^2}$.
The total number of samples can be bounded by 
\begin{align}\label{eq_sample_lspe_expectation}
\cO \rbr{
\frac{1}{\eta \mu} 
\rbr{1 + \frac{\gamma^2 \zeta^2 M^2 \overline{w}}{(1-\gamma)^2 \mu_w \epsilon^2}}
\log \rbr{\frac{r_\Theta }{ \epsilon}}
} .
% \cO \rbr{
%\frac{ L^2 + \kappa_2^2 + \kappa_3^2 + c_1^2}{ \mu^2} 
%\rbr{1 + \frac{\gamma^2 \zeta^2  M^2 \overline{w}}{(1-\gamma)^2 \mu \epsilon^2}}
%\log \rbr{\frac{r_\Theta }{(1-\gamma) \epsilon}}
%} .
\end{align}
In particular, when the distance generating function $w(\cdot)$ is set as in \eqref{dgf_negative_entropy}, the number of samples required is bounded by 
\begin{align*}
\cO \rbr{
\frac{1}{\eta \mu} 
\rbr{1 + \frac{\gamma^2 \zeta^2 M^2 \log \abs{\cS} }{(1-\gamma)^2  \epsilon^2}}
\log \rbr{\frac{r_\Theta }{\epsilon}}
} .
% \cO \rbr{
%\rbr{ L^2 + \kappa_2^2 + \kappa_3^2 + c_1^2}
%\rbr{1 + \frac{\gamma^2 \zeta^2  M^2   \log \abs{\cS} }{(1-\gamma)^2  \epsilon^2}}
%\log \rbr{\frac{r_\Theta }{(1-\gamma) \epsilon}}
%} .
\end{align*}
\end{theorem}

\begin{proof}
With a slight overload of notation, let us denote $\theta_k$ as the parameters output by SLPE at the $k$-th iteration of SFRPE. 
Then given the choice of parameters, one can apply Lemma \ref{lemma_ctd_recursion} and obtain  
%$
%\EE \sbr{ \norm{\theta_k - \theta^{\pi_k}}^2} \leq r_\Theta^2 + c_1^2
%$,
%and hence 
%\yan{ambiguity in notation $\theta_k$ here}
\begin{align*}
%\label{ctd_estimate_expectation_bound_on_korm}
\EE \sbr{ \norm{\cV^{\pi_k}_{\theta_k} - \cV^{\pi_k}}_\infty^2} 
\leq 2 \rbr{ \EE \norm{\theta_k - \theta^{\pi_k}}^2   + \varepsilon^2_{\mathrm{approx}}} 
\leq 2 \rbr{r_\Theta^2 + c_1^2  + \varepsilon^2_{\mathrm{approx}}}.
\end{align*}
Consequently, 
\begin{align*}
\EE  \sbr{ \norm{\cV^{\pi_k}_{\theta_k}}_\infty^2} \leq 2 \rbr{\EE \sbr{ \norm{\cV^{\pi_k}_{\theta_k} - \cV^{\pi_k}}_\infty^2} + \norm{\cV^{\pi_k}}_\infty^2 } \leq  4 \rbr{r_\Theta^2 + c_1^2  + \varepsilon^2_{\mathrm{approx}}} + {2}/{(1-\gamma)^2} .
\end{align*}
In addition, from Lemma \ref{lemma_bias_linear}, 
for any $\epsilon \geq \frac{8 \varepsilon_{\mathrm{approx}}}{1-\gamma} $, 
taking $T = \cO\rbr{ \frac{1}{\eta \mu} \log \rbr{\frac{r_\Theta }{ \epsilon}}}$ further yields that 
\eqref{stoch_expecation_conv_bias_condition} is satisfied with 
\begin{align*}
\frac{3\varepsilon}{1-\gamma} = \frac{3 \epsilon}{4}, ~ M^2 = 4 \rbr{r_\Theta^2 + c_1^2  + \varepsilon^2_{\mathrm{approx}}} + {2}/{(1-\gamma)^2} .
\end{align*}
Combining the above relation and Proposition \ref{thrm_stoch_generic_convergence_expectation}, the total number of iterations required by SFRPE is bounded by 
$
k = 1 + \frac{256 \gamma^2 \zeta^2 M^2 \overline{w}}{(1-\gamma)^2 \mu_w \epsilon^2}.
$
The total number of samples can be bounded by 
\begin{align*}
T \cdot k  & = \cO \rbr{
\frac{1}{\eta \mu} 
\rbr{1 + \frac{\gamma^2 \zeta^2 M^2 \overline{w}}{(1-\gamma)^2 \mu_w \epsilon^2}}
\log \rbr{\frac{r_\Theta }{\epsilon}}
} .
% \\
%& = \cO \rbr{
%\frac{ L^2 + \kappa_2^2 + \kappa_3^2 + c_1^2}{ \mu^2} 
%\rbr{1 + \frac{\gamma^2 \zeta^2  M^2 \overline{w}}{(1-\gamma)^2 \mu_w \epsilon^2}}
%\log \rbr{\frac{r_\Theta }{(1-\gamma) \epsilon}}
%} .
\end{align*}
The proof is then completed.
%We conclude the proof by noting  taking $\eta = \frac{\mu}{16(L^2 + \kappa_2^2 + \kappa_3^2 + c_1^2)}$ suffices to satisfy the condition \eqref{xx} of Lemma \ref{lemma_bias_linear}
\end{proof}

In view of Theorem \ref{thrm_lspe_expectation}, SFRPE yields an $\tilde{\cO}({\zeta^2}/ \epsilon^2 + \log(1/\epsilon))$ sample complexity with linear function approximation.
Notably, this appears to be the first method of stochastic robust policy evaluation beyond tabular settings with non-asymptotic convergence, and does so without restrictive assumptions on the transition kernel and discount factor \cite{tamar2014scaling, roy2017reinforcement}. 
Similar to Theorem \ref{thrm_sample_se_expectation} and \ref{thrm_stoch_se_high_prob}, the obtained sample complexity in \eqref{eq_sample_lspe_expectation} admits a natural decomposition, with the first term corresponding to learning the standard value function up to $\epsilon$-accuracy in bias, and the second term corresponding to the price of robustness.

With the same spirit as Theorem \ref{thrm_stoch_se_high_prob}, we proceed to show that with the same number of samples, SFRPE with SLPE operator can indeed output an $\epsilon$-estimator of the robust value function in high probability.
To this end, we first establish the following high probability bound on the iterate produced by SLPE operator. 
The challenge for such a statement primarily comes from that the noise in SLPE itself depends on the boundedness of the iterate, thus preventing a direct application of standard concentration argument. 
The following lemma constructs an approximate martingale sequence that upper bounds the noise of SLPE in high probability, with increment of the former sequence  bounded with proper choice of stepsize. 
Similar argument can also be found in \cite{li2023policy} for bounded increments but potentially non-zero conditional expectation.
%Such an argument can be viewed as a generalization of \cite{li2023policy} to unbounded sequence. 
%In view of this we proceed to show that the iterate of SLPE is indeed bounded in high probability. 
%In view of this observation, we adopt a similar technical idea that can be found in \cite{li2023policy}.
%, originally designed for establishing inherent exploration properties  of stochastic policy optimization.

\begin{lemma}\label{lemma_norm_bound_high_prob}
Fix total iterations $T > 0$.
For any $\delta \in (0,1)$, let $\eta_t = \eta \coloneqq \alpha / \sqrt{T}$ with 
\begin{align}
\alpha & \leq 
\min \big\{
%\frac{\mu \sqrt{T}}{\kappa_3^2}, ~
\frac{\mu }{32}, ~
\frac{1}{{192 \sbr{L^2 + c_1^2 + 16} }\log(2T / \delta)},~
\frac{1}{2 c_1}, ~
\frac{1}{4 G}
\big\},  \label{param_choice_bounded_norm_high_prob} \\
G & = 4 \sqrt{\log(2T /\delta)} \big[(L + 4) \norm{\theta_0 - \theta^{\pi}}^2 + c_1 \norm{\theta_0 - \theta^{\pi}}\big] + c_1. \label{G_bounded_norm_high_prob}
\end{align}
%where $G \geq 4 \sqrt{\log(2T /\delta)} \big[(L + \kappa_3) \norm{\theta_0 - \theta^{\pi}}^2 + (c_1 + \kappa_2  \varepsilon_{\mathrm{approx}}) \norm{\theta_0 - \theta^{\pi}}\big] + c_1$.
Then 
with probability at least $1-\delta$,
\begin{align*}
\norm{\theta_t - \theta^{\pi}}^2
\leq  \norm{\theta_0 - \theta^{\pi}}^2 + 1, ~ \forall t \leq T.
\end{align*}
\end{lemma}
%Going forward let us fix the epoch length $T > 0$ and the number of epochs $E > 0$.
%and assume $\eta_t = \eta \coloneqq \alpha / \sqrt{T}$ for some $\alpha > 0$ to be determined later. 
\begin{proof}
We begin by noting that  
\begin{align*}
\delta_t  = \hat{F}_t(\theta_t) - F(\theta_t)  =  \hat{F}_t(\theta_t)  - \hat{F}_t(\theta^\pi) -  F(\theta_t)  +  \hat{F}_t(\theta^\pi).
\end{align*}
Consequently, it is clear that  
\begin{align}\label{eq_norm_bound_with_noise_version}
\inner{\delta_t}{\theta_t - \theta^{\pi}}
& = 
 \inner{ \hat{F}_t(\theta_t)  - \hat{F}_t(\theta^\pi)}{\theta_t - \theta^\pi}
- \inner{F(\theta_t) }{\theta_t - \theta^\pi}
+ \inner{\hat{F}_t(\theta^\pi)}{\theta_t - \theta^\pi} \nonumber 
\\
& \leq 
(L + 4) \norm{\theta_t - \theta^\pi}^2  + 
 c_1 \norm{\theta_t - \theta^\pi},
\end{align}
where the last inequality applies \eqref{norm_bound_diff_stoch_op}, \eqref{norm_bound_stoch_op}, and \eqref{bound_on_op_norm}.
In addition, note that \eqref{norm_recursion_with_noise}  still holds.  Hence with 
$
 \eta \leq \frac{\mu}{32}, 
$
%or equivalently, $\alpha \leq \frac{\mu \sqrt{T}}{\kappa_3^2}$, 
or equivalently, $\alpha \leq \frac{\mu }{32}$, 
we obtain 
\begin{align}\label{recursion_on_norm_const_stepsize}
\norm{\theta_{t+1} - \theta^{\pi}}^2 
\leq  \norm{\theta_t - \theta^{\pi}}^2
+ 2\eta^2 c_1^2 +  2 \eta \inner{\delta_t}{\theta_t - \theta^{\pi}}.
\end{align}
Let us define random sequences $\cbr{X_t \coloneqq \inner{\delta_t}{\theta_t - \theta^{\pi}}}$, $\cbr{\tilde{X}_t \coloneqq \inner{\delta_t}{\theta_t - \theta^{\pi}} \mathbbm{1}_{\cG_t}}$,
where  $\cG_t = \cbr{Y_t \leq G \sqrt{t}}$, and
\begin{align*}
 Y_0 \equiv \tilde{Y}_0 \equiv 0,  ~
Y_t  = Y_{t-1} + X_{t-1}, ~
\tilde{Y}_t  = \tilde{Y}_{t-1} + \tilde{X}_{t-1}.
\end{align*}
%Now consider event $\cG_t = \cbr{Y_t \leq G \sqrt{te}}$ for some $G > 0$.
%Now consider event 
%$\cG_t = \cbr{\tsum_{i=0}^{t-1} \inner{\delta_i}{\theta_i - \theta^{\pi}} \leq G \sqrt{t}}$ for some $G > 0$.
Then over $ \cG_t$, 
taking the telescopic sum of \eqref{recursion_on_norm_const_stepsize}  yields   
\begin{align}
\norm{\theta_t - \theta^{\pi}}^2 & \leq \norm{\theta_0 - \theta^{\pi}}^2 + 2 t \eta^2 c_1^2 + 2 \eta G \sqrt{t} \nonumber \\
& \leq  \norm{\theta_0 - \theta^{\pi}}^2 + 2 \alpha^2 c_1^2  + 2 \alpha G  \coloneqq M_{(\alpha, G)}.
\label{norm_bound_on_good_event_generic}
\end{align}
We proceed to establish that \eqref{norm_bound_on_good_event_generic} indeed holds with probability at least $1-\delta$ given proper choice of $(\alpha, G)$.

%Let us define random sequences $\cbr{X_t \coloneqq \inner{\delta_t}{\theta_t - \theta^{\pi}}}$, $\cbr{\tilde{X}_t \coloneqq \inner{\delta_t}{\theta_t - \theta^{\pi}} \mathbbm{1}_{\cG_t}}$,
%together with 
%\begin{align*}
%Y_t = \tsum_{i=0}^{t-1} X_i, ~ \tilde{Y}_t = \tsum_{i=0}^{t-1} X_i', ~ Y_0 \equiv Y_0' \equiv 0.
%\end{align*}
%Then by definition we have $\cG_t = \cbr{Y_t \leq G \sqrt{t}}$.
First, it is  clear that 
\begin{align}\label{aux_sequence_norm_bound}
\abs{\tilde{X}_t} = \abs{\inner{\delta_t}{\theta_t - \theta^{\pi}} \mathbbm{1}_{\cG_t}} \leq (L + 4) M_{(\alpha, G)} + c_1 \sqrt{M_{(\alpha, G)}},
\end{align}
where the last inequality follows from \eqref{eq_norm_bound_with_noise_version} and \eqref{norm_bound_on_good_event_generic}.
Moreover, we also have 
\begin{align}
\abs{\EE_{|t}\sbr{\tilde{X}_t}}
& =  \abs{\EE_{|t}\sbr{\inner{\delta_t}{\theta_t - \theta^{\pi}}}  \mathbbm{1}_{\cG_t}}  = 0,\label{aux_sequence_conditional_expectation}
\end{align}
where the last equality follows from \eqref{bias_conditional_expectation}.
In view of \eqref{aux_sequence_norm_bound} and \eqref{aux_sequence_conditional_expectation}, we can now apply Azuma–Hoeffding inequality to $\cbr{\tilde{Y}_t}$ and obtain 
\begin{align}
\abs{
\tilde{Y}_t
} \leq 2 \sqrt{t}  \big[ (L + 4) M_{(\alpha, G)} + c_1 \sqrt{M_{(\alpha, G)}} \big]\sqrt{\log(2T/ \delta)}, 
~ \forall t \leq T,
\label{aux_sequence_accumulation_raw}
\end{align}
with probability $1-\delta$.
Through direct computation, one can verify that \eqref{aux_sequence_accumulation_raw} and the $(\alpha, G)$ specified in \eqref{param_choice_bounded_norm_high_prob}, \eqref{G_bounded_norm_high_prob}, together with the definition of $M_{(\alpha, G)}$, implies 
%\begin{align}
%& \sbr{ c_1 + \kappa_2  \varepsilon_{\mathrm{approx}} + L + \kappa_3} \sqrt{\alpha} \leq \frac{1}{8 \sqrt{\log(2T/\delta)}}, 
%~ \tau = \tilde{\cO} (t_{\mathrm{mix}}), \label{param_alpha_tau} \\
%& G \geq 4 \sqrt{\log(2T /\delta)} \big[(L + \kappa_3) \norm{\theta_0 - \theta^{\pi}}^2 + (c_1 + \kappa_2  \varepsilon_{\mathrm{approx}}) \norm{\theta_0 - \theta^{\pi}}\big] + c_1 , \label{param_choice_G}
%\end{align}
\begin{align}\label{aux_seq_accumulation_bounded}
\abs{\tilde{Y}_t} \leq G \sqrt{t }, ~ \forall t \leq T,
\end{align}
with probability $1-\delta$.
%and  $G = 4 \sqrt{\log(2k /\delta)} \sbr{(L + \kappa_3) \norm{\theta_0 - \theta^{\pi}}^2 + (c_1 + \kappa_2  \varepsilon_{\mathrm{approx}}) \norm{\theta_0 - \theta^{\pi}}} + c_1 + 1$.
Let us denote the event corresponding to \eqref{aux_seq_accumulation_bounded} by $\cG$.
%Clearly, from the definition of $\cbr{\cG_t}$ and \eqref{aux_seq_accumulation_bounded} we obtain $\cap_{t \leq T} \cG_t \subseteq \cG$.
Our next goal is to show that $Y_t = \tilde{Y}_t$ over $\cG$ for every $t \leq T$, and consequently 
\begin{align}\label{noise_accumulation_bound}
\abs{Y_t} \leq G \sqrt{t}, ~ \forall t \leq T, 
\end{align}
with probability $1-\delta$.
We proceed with an inductive argument. 
The claim trivially holds at $t = 0$.
Suppose the claim holds at iteration $t \geq 0$, then for any $\omega \in \cG$, 
\begin{align}
Y_{t+1}(\omega) & = Y_t(\omega) + X_t(\omega) \nonumber \\
& = Y_t(\omega) + X_t(\omega) \mathbbm{1}_{\cG}(\omega) \nonumber \\
& \overset{(a)}{=} Y_t(\omega) + X_t(\omega) \mathbbm{1}_{\cbr{\tilde{Y}_t \leq G \sqrt{t}}}(\omega) \nonumber \\
& \overset{(b)}{=} Y_t(\omega) + X_t(\omega) \mathbbm{1}_{\cbr{Y_t \leq G \sqrt{t}}}(\omega) \nonumber \\
& \overset{(c)}{=} Y_t(\omega) + \tilde{X}_t(\omega) \nonumber \\
& = \tilde{Y}_{t+1}(\omega), \label{equiv_induction_step_1}
\end{align}
where $(a)$ follows from the definition of $\cG$ and \eqref{aux_seq_accumulation_bounded}, 
 $(b)$ follows from the induction hypothesis that $Y_t (\omega) = \tilde{Y}_t(\omega)$ for $\omega \in \cG$,
and $(c)$ applies the definition of $\tilde{X}_t$.
Hence \eqref{equiv_induction_step_1} completes the induction step.
%On the other hand, suppose the claim holds for all $t \leq T$ of epoch $e$, then for any $\omega \in \cG$,
%\begin{align*}
%Y_{0}^{(e+1)}(\omega) & = Y_{T}(\omega) = \tilde{Y}_{T}(\omega)  = \tilde{Y}^{(e+1)}_{0}(\omega).
%\end{align*}
%\begin{align*}
%Y_{0}^{(e+1)}(\omega) & = Y_{T}(\omega) \\
%& = Y_{T-1}(\omega) + X_{T-1}(\omega) \mathbbm{1}_{\cG}(\omega) \\
%& = Y_{T-1}(\omega) + X_{T-1}(\omega) \mathbbm{1}_{\cbr{\tilde{Y}_{T-1} \leq G \sqrt{T-1 + e T}}}(\omega) \\
%& =Y_{T-1}(\omega) + X_{T-1}(\omega) \mathbbm{1}_{\cbr{Y_{T-1} \leq G \sqrt{T-1 + e T}}}(\omega) \\
%& = Y_{T-1}(\omega) + \tilde{X}_{T-1}(\omega) \\
%& = \tilde{Y}^{(e+1)}_{0}(\omega),
%\end{align*}
%Combining \eqref{equiv_induction_step_1} and the above relation completes the induction step.

In view of \eqref{norm_bound_on_good_event_generic} and \eqref{noise_accumulation_bound}, we conclude that
% for $(\alpha, \tau, G)$ specified  in  \eqref{param_alpha_tau} and \eqref{param_choice_G}, 
\begin{align*}
\norm{\theta_t - \theta^{\pi}}^2
\leq \norm{\theta_0 - \theta^{\pi}}^2 + 2 \alpha^2 c_1^2  + 2 \alpha G,
\end{align*}
with probability at least $1-\delta$.
The desired claim follows immediately by noting that $\alpha \leq \min\cbr{\frac{1}{2 c_1}, ~
\frac{1}{4 G} }$.
\end{proof}
%With \eqref{xx}, \eqref{xx} and \eqref{xx} in place, we are now ready to complete the proof.
%Clearly, with an inductive argument, one can readily show that with 
%\begin{align*}
%& \sbr{ c_1 + \kappa_2  \varepsilon_{\mathrm{approx}} + L + \kappa_3} \sqrt{\alpha} \leq \frac{1}{8 \sqrt{\log(2T/\delta)}}, ~
%\alpha c_1 \leq 1, 
%~ \tau = \cO(t_{\mathrm{mix}}),  \\
%& G \geq 4 \sqrt{\log(2k /\delta)} \big[(L + \kappa_3) (e + 1) \norm{\theta_0 - \theta^{\pi}}^2 + (c_1 + \kappa_2  \varepsilon_{\mathrm{approx}}) \sqrt{e + 1} \norm{\theta_0 - \theta^{\pi}}\big] + c_1 + 1, \\
%& 
% 2 \alpha^2 c_1^2  + 2 \alpha G \leq \norm{\theta_0 - \theta^{\pi}}^2.
%\end{align*}

By combining Lemma \ref{lemma_bias_linear} and \ref{lemma_norm_bound_high_prob}, we proceed to establish that SLPE with proper parameter specification yields fast bias reduction while controlling boundedness of the estimated value function $\hat{\cV}^\pi$.

\begin{lemma}\label{lemma_high_prob_norm_and_bias}
Fix total iterations $T > 0$ a priori in SLPE.
For any $\delta \in (0,1)$ and any $\varepsilon \geq 2 \varepsilon_{\mathrm{approx}}$, let the parameters in SLPE be chosen as 
\begin{align*}
\eta_t = \alpha /\sqrt{T}, 
%~ T = \frac{2 \log^2 \rbr{4{\norm{\theta_0 -\theta^{\pi}}}/{\varepsilon}}}{\alpha^2 \mu^2}, 
\end{align*} 
with 
$
\alpha  \leq 
\min \big\{
\frac{\mu }{L^2 + 32}, 
%\frac{\mu \sqrt{T}}{\kappa_3^2}, ~
\frac{1}{192 \sbr{L^2 + c_1^2 + 16} \log(2T / \delta)},
\frac{1}{2 c_1}, 
\frac{1}{4 G}
\big\},
%\\
%G & \geq 4 \sqrt{\log(2T /\delta)} \big[(L + \kappa_3) \norm{\theta_0 - \theta^{\pi}}^2 + (c_1 + \kappa_2  \varepsilon_{\mathrm{approx}}) \norm{\theta_0 - \theta^{\pi}}\big] + c_1. 
$
and $G$ defined in \eqref{G_bounded_norm_high_prob}.
Then the number of iterations required by SLPE to output 
\begin{align}\label{ctd_bias_bd_with_norm_bd}
\norm{\EE \sbr{ \cV^{\pi}_{\theta_t}} - \cV^{\pi} }_\infty \leq  \varepsilon
\end{align}
  is bounded by
\begin{align}\label{ctd_num_iter_bias_and_norm}
T = \cO \rbr{
 \frac{\log^2 \rbr{{ \norm{\theta_0 -\theta^{\pi}}}/{\varepsilon}}}{\alpha^2 \mu^2} 
 }.
%\end{align}
%The number of samples can be bounded by 
%\begin{align}
%\label{ctd_num_samples_bias_and_norm}
% \cO \rbr{
%\frac{t_{\mathrm{mix}}}{\alpha^2 \mu^2} \log^2 \rbr{\frac{\norm{\theta_0 -\theta^{\pi}}}{ \varepsilon}}
%}.
\end{align}
In addition, we have 
\begin{align}\label{ctd_norm_bound_bias_and_norm}
\norm{\cV^{\pi}_{\theta_t} -  \cV^{\pi} }_\infty
\leq  \norm{\theta_0 - \theta^{\pi}} + 1 + \varepsilon_{\mathrm{approx}}, ~ \forall t \leq T ,
\end{align}
with probability at least $1-\delta$.
\end{lemma}

\begin{proof}
Clearly,  $\eta = \alpha/\sqrt{T}$ with the choice of specified $\alpha$ satisfy 
 \eqref{param_choice_bias_reduction}.
 Consequently, one can  apply \eqref{bound_on_value_bias} in Lemma \ref{lemma_bias_linear}, and obtain that for any $\varepsilon \geq 2 \varepsilon_{\mathrm{approx}}$, 
SLPE outputs 
$
\norm{\EE \sbr{ \cV^{\pi}_{\theta_t}} - \cV^{\pi}}_\infty \leq \varepsilon
$
in 
\begin{align*}
T =\cO \rbr{\frac{1}{\eta \mu} \log \rbr{\frac{\norm{\theta_0 - \theta^*}}{\varepsilon}}}
\end{align*} 
steps. Combining the above relation with the definition of $\eta = \alpha /\sqrt{T}$ implies \eqref{ctd_num_iter_bias_and_norm}.
%\eqref{ctd_num_samples_bias_and_norm} then follows from \eqref{ctd_num_iter_bias_and_norm} and the choice of $\tau = \tilde{\cO}(t_{\mathrm{mix}})$. 
%The total number of samples consumed by SLPE is bounded by 
%\begin{align*}
%T \cdot \tau 
%= \cO \rbr{
%\frac{t_{\mathrm{mix}}}{\alpha^2 \mu^2} \log^2 \rbr{\frac{\norm{\theta_0 -\theta^{\pi}}}{\varepsilon}}
%}.
%\end{align*}
Finally, \eqref{ctd_norm_bound_bias_and_norm} follows  from 
 \begin{align*}
  |\cV^{\pi}_{\theta_t} (s) - \cV^{\pi}  |_\infty \leq \norm{\theta_t - \theta^*} +  \varepsilon_{\mathrm{approx}} \leq  \norm{\theta_0 - \theta^{\pi}} + 1 + \varepsilon_{\mathrm{approx}}, ~ \forall s \in \cS,
 \end{align*}
 where the last inequality applies Lemma \ref{lemma_norm_bound_high_prob}.
 The proof is then completed.
\end{proof}

%
%\begin{remark}
%One can also readily employ Lemma \ref{lemma_ctd_recursion} and \ref{lemma_norm_bound_high_prob} to establish the convergence of $\EE \sbr{\norm{\overline{\theta}_t -\theta^{\pi}}^2} $ with high probability, for proper defined ergodic iterate $\overline{\theta}_t$.
%We omit its explicit discussion to keep the scope of the manuscript concise. 
%\yan{remark on extension to least square to improve the dependence on visitation measure, with simplified analysis. also this can exploits the generator}
%\end{remark}

With Lemma \ref{lemma_bias_linear} and \ref{lemma_norm_bound_high_prob} in place, we can now establish the sample complexity of SFRPE with the SLPE operator 
that outputs an $\epsilon$-estimator of the robust value function
with high probability.

\begin{theorem}\label{thrm_sample_slpe_high_prob}
Fix total iterations $k > 0$ in SFRPE and $\delta \in (0,1)$.
For any $\epsilon \geq \frac{8 \varepsilon_{\mathrm{approx}} }{ (1-\gamma)}$, 
 let SFRPE be instantiated with the SLPE operator with evaluation parameters 
\begin{align*}
\eta_t = \alpha /\sqrt{T},  ~
\theta_0 = 0, 
~ T = \cO \rbr{
 \frac{\log^2 \rbr{{ r_\Theta}/{\epsilon}}}{\alpha^2 \mu^2} 
 }, 
\end{align*} 
where 
\begin{align*}
\alpha  & \leq 
\min \big\{
\frac{\mu }{L^2 + 32}, ~
%\frac{\mu \sqrt{T}}{\kappa_3^2}, ~
\frac{1}{192 \sbr{L^2 + c_1^2 + 16} \log(12 Tk/ \delta)},~
\frac{1}{2 c_1}, ~
\frac{1}{4 G}
\big\}
\\
G&   = 4 \sqrt{\log(12 Tk/\delta)} \big[(L + 4) r_{\Theta}^2 + c_1   r_\Theta \big] + c_1,
\end{align*}
and optimization parameters specified as
\begin{align*}
\beta_t = t^{1/2},   ~ \lambda_t =  \frac{ (t+1) M \gamma \zeta }{2  \sqrt{\mu \overline{w}}}, ~\forall t \leq k,
\end{align*}
where $M \geq r_\Theta + 1 + \varepsilon_{\mathrm{approx}} + \frac{1}{1-\gamma}$.
Then with probability at least $1-\delta$ we have 
\begin{align}
& -\frac{\epsilon}{4} -  \frac{4M}{\sqrt{k}} \sqrt{\log (\frac{6 (k+1)\abs{\cS}}{\delta})} \\
  \leq &  \tsum_{t=1}^k \theta_t 
\cV^{\pi_t}(s) - \cV^*(s) \nonumber \\
 \leq & \frac{4 \gamma \zeta M \sqrt{\overline{w}}}{(1-\gamma) \sqrt{\mu k}}  + \frac{3 \epsilon }{4}
+ \frac{8 \gamma \zeta M}{(1-\gamma) \sqrt{k}} \sqrt{\log (\frac{6 (k+1)\abs{\cS}}{\delta})}
+ \frac{4M}{\sqrt{k}} \sqrt{\log (\frac{6(k+1) \abs{\cS}}{\delta})}, ~ \forall s \in \cS,
 \label{high_prob_err_bound_ctd}
\end{align}
The total number of samples required by SFRPE to output $-  \epsilon \leq \tsum_{t=1}^k \theta_t
\cV^{\pi_t}(s) - \cV^{\pi^*}(s) \leq  \epsilon$ for all $s \in \cS$ with at least probability $1-\delta$  is bounded by 
\begin{align}\label{ctd_sample_high_prob}
\tilde{\cO} \rbr{
\frac{1}{\alpha^2 \mu^2}
\rbr{
\frac{\gamma^2 \zeta^2 M^2 \overline{w}}{(1-\gamma)^2 \mu_w \epsilon^2}
+ \frac{M^2}{\epsilon^2}
}
 \log^2 \rbr{\frac{r_\Theta}{\epsilon}}
}.
\end{align}
In particular, when the distance generating function $w_s(\cdot)$ is set as in \eqref{dgf_negative_entropy}, the total number of samples required can be bounded by 
\begin{align*}
\tilde{\cO} \rbr{
\frac{1}{\alpha^2 \mu^2}
\rbr{
\frac{\gamma^2 \zeta^2 M^2 \log \abs{\cS}}{(1-\gamma)^2  \epsilon^2}
+ \frac{M^2}{\epsilon^2}
}
 \log^2 \rbr{\frac{r_\Theta}{\epsilon}}
}.
\end{align*}
\end{theorem}

\begin{proof}
The essential argument is similar to that of Theorem \ref{thrm_stoch_se_high_prob}, but we will use Lemma \ref{lemma_high_prob_norm_and_bias} instead of Proposition \ref{thrm_stoch_generic_convergence_expectation}.
%Denote $M = r_\Theta + 1 + \varepsilon_{\mathrm{approx}}$.
Clearly the choice of parameters satisfies conditions of Lemma \ref{lemma_high_prob_norm_and_bias}, and hence
\begin{align}\label{norm_bound_every_iter_whp}
\norm{\cV^{\pi_t}_{\theta_t}  }_\infty \leq M, ~
\norm{\EE \sbr{ \cV^{\pi_t}_{\theta_t}} - \cV^{\pi_t}}_\infty \leq   \frac{(1-\gamma) \epsilon}{4},  ~ \forall t \leq k, 
\end{align}
 with probability $1 -  \delta / 3$.
Combining the above relation with Proposition 34 of \cite{tao2015random} yields
\begin{align}\label{ctd_noise_accumulation_1}
\abs{ \tsum_{t=1}^k \theta_t \rbr{\cV^{\pi_t}(s) - \cV^{\pi_t}(s)}  } 
  \leq 
\frac{(1-\gamma) \epsilon}{4} 
+  M \sqrt{2 \tsum_{t=1}^k \theta_t^2 \log(\frac{2}{\delta})}  \leq \frac{(1-\gamma) \epsilon}{4} + \frac{4 M }{\sqrt{k}} \sqrt{\log(\frac{2 }{\delta})}, 
\end{align}
with probability at least $1 - \delta / 6$, for every $s\in \cS$. Further applying union bound yields 
\begin{align*}
\abs{ \tsum_{t=1}^k \theta_t \rbr{\cV^{\pi_t}(s) - \cV^{\pi_t}(s)}  } 
  \leq \frac{(1-\gamma) \epsilon}{4} + \frac{4 M }{\sqrt{k}} \sqrt{\log(\frac{2 \abs{\cS} }{\delta})}, ~ \forall s \in \cS, 
\end{align*}
with probability at least $1 -  \delta/3$.
Applying the same treatment, one can also show that 
\begin{align}\label{ctd_noise_accumulation_2}
\tsum_{t=1}^k \frac{\theta_t}{1-\gamma}  \EE_{s' \sim d_{s}^{\pi^*}} \sbr{ \delta_t(s', \pi^*(s'))}
 & \leq 
 \frac{\gamma \zeta}{1-\gamma} 
\rbr{
 \frac{(1-\gamma) \epsilon}{2}
+ \frac{8 M }{ \sqrt{k}} \sqrt{  \log(\frac{2 \abs{\cS} }{\delta})}
}, ~ \forall s \in \cS,
\end{align}
with probability at least $1 - \delta/3$. 
By plugging \eqref{norm_bound_every_iter_whp},  \eqref{ctd_noise_accumulation_1} and \eqref{ctd_noise_accumulation_2} into Lemma \ref{lemma_generic_prop_stoch}, we obtain
\begin{align*}
-\frac{\epsilon}{4} -  \frac{4M}{\sqrt{k}} \sqrt{\log (\frac{\abs{2 \cS}}{\delta})}
&  \leq \tsum_{t=1}^k \theta_t 
\cV^{\pi_t}(s) - \cV^{\pi^*}(s) \nonumber \\
& \leq \rbr{\tsum_{t=1}^k \beta_t}^{-1}  \tsum_{t=1}^k \frac{\beta_t^2 \gamma^2 \zeta^2  M^2}{2 \mu \lambda_{t-1} (1-\gamma)}
+  \rbr{\tsum_{t=1}^k \beta_t}^{-1} \frac{\lambda_k \overline{w}}{1-\gamma} \nonumber \\
& ~~~ +  \frac{3 \epsilon }{4}
+ \frac{8 \gamma \zeta M}{(1-\gamma) \sqrt{k}} \sqrt{\log (\frac{2 \abs{\cS}}{\delta})}
+ \frac{4M}{\sqrt{k}} \sqrt{\log (\frac{2 \abs{\cS}}{\delta})}, ~ \forall s \in \cS,
\end{align*}
with probability $1-  \delta $. 
Plugging the choice of $\cbr{(\beta_t, \lambda_t)}$ into the above relation yields \eqref{high_prob_err_bound_ctd}.
%\begin{align*}
%& -\epsilon -  \frac{4M}{\sqrt{k}} \sqrt{\log (\frac{6 (k+1)\abs{\cS}}{\delta})} \\
%  \leq &  \tsum_{t=1}^k \theta_t 
%\cV^{\pi_t}(s) - \cV^{\pi^*}(s) \nonumber \\
% \leq & \frac{4 \gamma \zeta M \sqrt{\overline{w}}}{(1-\gamma) \sqrt{\mu k}}  + \rbr{1 + \frac{2}{1-\gamma}} \epsilon 
%+ \frac{8 \gamma \zeta M}{(1-\gamma) \sqrt{k}} \sqrt{\log (\frac{6 (k+1)\abs{\cS}}{\delta})}
%+ \frac{4M}{\sqrt{k}} \sqrt{\log (\frac{6(k+1) \abs{\cS}}{\delta})}, 
%\end{align*}
%for any $s \in \cS$, with probability at least $1-\delta$.
Given \eqref{high_prob_err_bound_ctd}, the total number of iterations by SFRPE  to output $-  \epsilon \leq \tsum_{t=1}^k \theta_t 
\cV^{\pi_t}(s) - \cV^{\pi^*}(s) \leq  \epsilon$ can be bounded by 
\begin{align*}
k = \tilde{\cO} \rbr{
\frac{\gamma^2 \zeta^2 M^2 \overline{w}}{(1-\gamma)^2 \mu_w \epsilon^2}
+ \frac{M^2}{\epsilon^2}
}.
\end{align*}
We conclude the proof by noting that the total number of samples required is $k T$.
%The total number of samples is bounded by 
%\begin{align*}
%\cO \rbr{
%\frac{t_{\mathrm{mix}}}{\alpha^2 \mu^2}
%\rbr{
%\frac{\gamma^2 \zeta^2 M^2 \overline{w}}{(1-\gamma)^2 \mu_w \epsilon^2}
%+ \frac{M^2}{\epsilon^2}
%}
% \log^2 \rbr{\frac{r_\Theta}{\epsilon}}
%}
%\end{align*}
%Finally, \eqref{ctd_sample_high_prob} follows from the above relation and \eqref{ctd_num_samples_bias_and_norm}.
\end{proof}

In view of Theorem \ref{thrm_sample_slpe_high_prob}, SFRPE method instantiated with LSPE operator attains an $\tilde{\cO}(\zeta^2/\epsilon^2 + 1/\epsilon^2)$ sample complexity to output an $\epsilon$-estimator of the robust value function in high probability. 
Clearly the accuracy certificate of Theorem \ref{thrm_sample_slpe_high_prob} is stated in a stronger sense compared to the expectation statement in Theorem \ref{thrm_lspe_expectation}.
As before, the sample complexity possesses two terms that can be attributed to the price of robustness and standard value function estimation, respectively.

\vspace{6pt}
%Before we conclude our discussion in this section.
%\begin{remark}[Applications to Stochastic Policy Optimization for Robust MDPs]
{\bf  Applications to Stochastic Policy Optimization for Robust MDPs.}
We conclude this section by briefly demonstrating the application of the developed results in the context of  stochastic policy optimization for robust MDPs. 
Consider solving large-scale robust MDPs with $(\mathrm{s}, \mathrm{a})$-rectangular sets using the stochastic robust policy mirror descent (SRPMD) method in \cite{li2022first}.
If log-linear policy class is employed, and one applies the same feature map for the policy class and the robust state-action value function, then applying SFRPE for learning the robust state-action value function immediately implies an $\tilde{\cO}(1/\epsilon^2)$ sample complexity of SRPMD for finding an $\epsilon$-optimal robust policy (cf. Proposition 5.1, \cite{li2022first}).
In particular, this simple application already yields the first sample complexity of policy gradient methods applied to robust MDPs beyond tabular~settings.

\section{Concluding Remarks}\label{sec_conclusion}

We introduce a first-order method, named FRPE, applied to robust policy evaluation problem by taking a policy optimization viewpoint.
For offline robust MDPs, FRPE attains linear convergence, and is compatible with function approximation for large-scale MDPs. 
For online robust MDPs, we further develop a stochastic variant that attains optimal $\cO(1/\epsilon^2)$ sample complexity.
The obtained sample complexities also delineate the clear role of robustness when performing policy evaluation. 
Notably these sample complexities are established for both tabular setting and for large-scale problems where linear approximation is applied. 
We conclude by discussing the application of FRPE in the context of stochastic policy optimization applied to robust MDPs. 
It appears to us that the FRPE framework overcomes long-standing challenges of robust policy evaluation and opens up a new spectrum of potential directions. 
We now discuss a few of them as follows. 

First, though we mainly presents linear function approximation for stochastic robust policy evaluation,  it should be noted stochastic FRPE can be combined with generic function approximation scheme as long as the stochastic evaluation of $\cV^\pi$ can be performed efficiently.
As an example one can consider employing over-parameterized neural networks in SLPE. With recent advances in showing their global convergence, it appears that a polynomial sample complexity can be attained. 

Second, the proposed LSPE for estimating $\cV^\pi$ requires the simulator access. This seems to be natural for online robust MDPs motivated by sim-to-real problems (cf. Example \ref{online_robust_mdp}), where a simulator is inherently available. 
As the estimation of $\cV^\pi$ can be essentially viewed as an off-policy evaluation problem (cf., \eqref{eq_nature_value_as_player_value}), one can also apply existing off-policy evaluation methods that can be implemented without simulator access while being compatible with linear function approximation \cite{sutton2008convergent, sutton2009fast}. 
To attain an $\tilde{\cO}(1/\epsilon^2)$ sample complexity for stochastic FRPE, it is important to develop two properties of the method. 
Namely the fast bias reduction (Lemma \ref{lemma_bias_linear}) and the high probability control on the boundedness of the iterate (Lemma \ref{lemma_norm_bound_high_prob}). 
%For this purpose we believe a single time-scale modification is required before adopting the aforementioned methods that are two-timescale in nature.

Third, while  $\cO(1/\epsilon^2)$ samples suffices for stochastic FRPE  to output an $\epsilon$-estimator of the robust value function in expectation, it remains to be seen whether this sample complexity is tight. 
In particular, when $\zeta = 0$, and the robust policy evaluation reduces to standard policy evaluation,  it can be seen from Theorem \ref{thrm_sample_se_expectation} and \ref{thrm_lspe_expectation} that only $\cO(\log(1/\epsilon))$ samples are needed. 
Whether the phase transition between robust and standard policy evaluation exists could be impactful for potentially improving the sample complexity of stochastic policy optimization applied to robust MDPs.
We refer interested readers to \cite{li2022first} for related~discussions.

Finally, it is also highly rewarding to extend FRPE to a broader class of ambiguity sets  beyond $\mathrm{s}$- and $(\mathrm{s}, \mathrm{a})$-rectangularity \cite{mannor2016robust, goyal2023robust},
and to settings where ambiguity sets admit more complex structure than \eqref{def_ambiguity_set_structure}. 
%and discuss their applications in the context of policy optimization.

%\yan{
%Potential directions:
%\begin{itemize}
%\item the framework allows using more advanced function approximation -- neural networks (Zhaoran's paper)
%\item purely online evaluation method for LSPE, for example, GTD 
%\item whether bias can be shown to converge faster 
%\item application of the results in policy optimization methods
%\item other ambiguity sets
%\end{itemize}
%
%}

%\newpage
\bibliographystyle{plain}
\bibliography{references}

\end{document}